\documentclass{amsart}
\usepackage{lineno,hyperref}
    \hypersetup{
        colorlinks   = true,
        citecolor    = blue
    }
\usepackage[utf8]{inputenc}
\usepackage{subcaption}
\usepackage{algorithmicx}
\usepackage{algorithm}
\usepackage{algpseudocode}
\usepackage{soul} 
\usepackage{url}
\usepackage[left=2.5cm,top=2.5cm,right=2.5cm,bottom=2.5cm]{geometry}

\usepackage{enumerate}
\usepackage{lipsum}
\usepackage{amsfonts}
\usepackage{latexsym, color}
\usepackage{colordvi}
\usepackage{amsmath}
\usepackage{longtable}
\usepackage{multirow}
\usepackage{enumitem}
\usepackage{tabularx}
\usepackage{tabulary}
\usepackage{natbib}
\bibpunct{(}{)}{;}{a}{,}{;}

\usepackage{booktabs}
\usepackage{adjustbox}

\begin{document}

\title[GRASP with PR - 20 years]{20 years of Greedy Randomized Adaptive Search\\ Procedures with Path Relinking}

\author[M. Laguna]{M. Laguna}
\address[Manuel Laguna]{
Leeds School of Business, University of Colorado Boulder (USA)
}
\email[]{laguna@colorado.edu}

\author[R. Mart\'{\i}]{R. Mart\'{\i}}
\address[Rafael Mart\'{\i}]{Departament d'Estad\'{\i}stica i Investigaci\'o Operativa,  Universitat de Val\`encia (Spain)}
\email[]{rafael.marti@uv.es}

\author[A. Martinez-Gavara]{A. Martinez-Gavara}
\address[Anna Martinez-Gavara]{Departament d'Estad\'{\i}stica i Investigaci\'o Operativa,  Universitat de Val\`encia (Spain)}
\email[]{gavara@uv.es}

\author[S. Perez-Pel\'o]{S. Perez-Pel\'o}
\address[Sergio Perez-Pel\'o]{Department of Computer Science, 
University Rey Juan Carlos (Spain)}
\email[]{sergio.perez.pelo@urjc.es}

\author[M.G.C. Resende]{M.G.C. Resende}
\address[Mauricio G.C. Resende, \textcolor{black}{ORCID 0000-0001-7462-6207}]{Industrial \& Systems Engineering, University of Washington, 
Seattle, WA (USA)}
\email[M.G.C. Resende]{mgcr@uw.edu}

\begin{abstract}
This is a comprehensive review of the Greedy Randomized Adaptive Search Procedure (GRASP) metaheuristic and its hybridization with Path Relinking (PR) over the past two decades. GRASP with PR has become a widely adopted approach for solving hard optimization problems since its proposal in 1999. The paper covers the historical development of GRASP with PR and its theoretical foundations, as well as recent advances in its implementation and application. The review includes a critical analysis of variants of PR, including memory-based and randomized designs, with a total of ten different implementations. It describes these advanced designs both theoretically and practically on two well-known optimization problems, linear ordering and max-cut. The paper also explores the hybridization of GRASP with PR and other metaheuristics, such as Tabu Search and Scatter Search. Overall, this review provides valuable insights for researchers and practitioners seeking to utilize GRASP with PR for solving optimization problems.
\end{abstract}
\date{December 18, 2023}
\maketitle

\section{The origins of the methodology}\label{sec:origins}

The GRASP methodology was developed in the late 1980s \citep{feo1989}, and stands for Greedy Randomized Adaptive Search Procedure \citep{res2016}. During each iteration of GRASP, a trial solution is constructed, and local search is then applied to improve it. The construction phase is iterative, pseudo-greedy, and adaptive. GRASP independently samples the solution space, and no information, other than the incumbent solution, is passed from one iteration to the next.

Path relinking (PR) was initially proposed as a mechanism for long-term memory within tabu search. In GRASP, PR enhances the memory component by exploring trajectories that link new GRASP solutions to elite solutions discovered during the search. This strategy intensifies the search by generating new solutions through paths in the neighborhood space, starting with an initiating solution and moving towards the guiding solutions. PR selects moves to introduce attributes from the guiding solutions. In tabu search, any two solutions are linked by a  series of moves executed during the search. The path between solutions is determined by the ``normal'' operational rules of the search procedure. For example, neighborhood searches of the type typically implemented within tabu search transition from one solution to another by greedy rules that seek the maximum objective function improvement in the local sense. The relinking process, on the other hand, seeks to create a new path by incorporating attributes of a guiding (elite) solution. The moves chosen during the relinking process are different from those moves during the ``normal'' search because the relinking moves do not use the change of the objective function as the guiding principle.

The idea of combining GRASP and path relining originated in the context of a 2-layer straight line crossing minimization problem \citep{lag1999}. This problem arises in the context of automated drawing systems. Developers of these systems are concerned with drawing speed and therefore the research efforts in line crossing minimization algorithms have focused on simple heuristic procedures that sacrifice solution quality for speed. The goal in this seminal paper was to develop a procedure that could compete on speed with the simple heuristics and on solution quality with metaheuristics whose solution time was impractical for automated drawing systems. To achieve this balance between solution quality and speed, the GRASP implementation included a mechanism that estimated the likelihood that the local search could produced an improvement when starting from a newly constructed solution. This filtering of solutions to be subjected to the local search was an important element of the implementation because path relinking is applied to all local optima. The ``relinking" concept has a special meaning within GRASP because local optima found in separate GRASP iterations are not linked by a sequence of moves (like in tabu search). The method implemented in \citet{lag1999} is the so-called forward path relinking in which a newly found local optimum is used as the initiating solution and an elite solution (chosen randomly from a set of the best three solutions found so far) as the guiding solution. At each step of the process, the attribute from the guiding solution (not present in the initiating solution) that results in the largest improvement is incorporated into the initiating solution. Throughout this paper, we consider maximization problems. In those problems, such as the ones presented in subsequent sections of this paper, improvement is measured as an increase in the objective function value.

\subsection{GRASP}

Many metaheuristic algorithms rely on local search for intensification.  Local search is a systematic search in a local neighborhood of a given solution $S$.  A local neighborhood $N(S)$ of $S$ is defined as the set of solutions that can be obtained from $S$ by making a small perturbation in $S$, often referred to as a \textit{move}.  For example, if $S$ is represented as a permutation, then one neighborhood of $S$ may be the set of solutions obtained by exchanging two elements of the permutation. If $S$ is represented as a 0-1 indicator vector, a neighborhood of $S$ may be the set of all indicator vectors obtained by flipping the value of one of its elements.

Local search starts from an initial solution $S_0$ and explores the solutions in $N(S_0)$ seeking a solution with better objective value.  In a first-improving variant, local search resets $S_0$ to be the first improving solution $S_*$ found in $N(S_0)$, i.e. $S_0 \gets S_*$ and the local search is restarted at $S_*$.   In a best-improving variant, local search resets $S_0$ to be the best improving solution found in $N(S_0)$ and restarts. If no improving solution is found in $N(S_0)$, then we say that $S_0$ is a locally optimal solution and the local search halts.  

By applying local search from different starting solutions $S_0$ a variety of locally optimal solutions may be found.  Embedding local search within a multi-start procedure where each local search starts from a different starting solution will produce a set of locally optimal solutions, the best of which could perhaps be a global optimum.

GRASP is based on the idea of generating a new starting solution for local search by using a randomized greedy, or semi-greedy, algorithm to generate $S_0$.  A greedy algorithm builds a solution one element at a time. At each step of the greedy construction a candidate element is chosen and added to the solution.  Candidate elements are those which if added to the partially constructed solution either reduce or do not increase infeasibility.
Among all candidate elements, the greedy algorithm selects an element which has the maximum contribution to the objective function value. Once the greedy candidate has been added to the partial solution, then one or more of the remaining candidate elements may have to be removed from the candidate set if their inclusion in the solution increases infeasibility. The contributions of the remaining candidate elements may need to be updated to take into account the newly added element. Applying the greedy algorithm repeatedly will always lead to the same greedy solution (assuming ties are broken deterministically) and local search would also lead a single solution.
The GRASP approach is to randomize the greedy construction to obtain multiple local optima.  Instead of selecting an element with the maximum contribution to the partial solution, GRASP builds a set of candidates with high contribution and selects an element from this set at random. This set is referred to in \citet{FeoRes1995a} as the Restricted Candidate List (RCL) and it can be value based or cardinality based as described below.

A value-based RCL involves setting a threshold based on the quality or value of the candidate elements. Only those candidates whose contribution to the solution under construction is within a certain range of the best-valued choice are included in the RCL. For instance, the threshold can be set as a percentage of the greedy choice or a fixed value difference from the greedy choice. The candidates that meet this criterion are then considered for random selection. This threshold is usually determined by a parameter $\alpha$, where $0 \leq \alpha \leq 1$.  By varying the value of $\alpha$ one can control the amount of randomness or greediness in the construction.

In a cardinality-based approach, the RCL is defined by a fixed number of top candidates. Instead of using a value-based threshold, the RCL consists of a predetermined number $k$ of the best candidates. The size of the list is constant, and candidates are chosen based on their ranking in terms of solution quality. The top candidates (e.g., top 5, top 10) are included in the RCL for random selection. The amount of randomness or greediness within the construction in this approach is also controlled by a parameter, in this case $k$.
Both approaches aim to strike a balance between exploration and exploitation in the search space. Larger restricted candidate lists emphasize exploration while smaller lists favor exploitation.

GRASP is a multi-start metaheuristic algorithm where, as shown in Algorithm~\ref{alg:alg-1} for the case of maximization, each iteration consists of a greedy randomized construction followed by local search.  The best locally optimal solution is returned as the GRASP solution.

\begin{algorithm}[htbp]
	\caption{\texttt{GRASP}($\alpha $)}
	\label{alg:alg-1}
	\begin{algorithmic}[1]
           \State{$f_* \gets -\infty$} \label{alg:grasp:sb}
            \While {\textrm{stopping criterion is not satisfied}} \label{alg:grasp:forini}
                \State{$ S \gets \texttt{Greedy\_Randomized\_Adaptive\_Construction}(\alpha) $} \label{alg:grasp:grc}
                \State{$ S \gets \texttt{Local\_Search}(S) $} \label{alg:grasp:li}
                \If{$f(S) > f_*$} \label{alg:grasp:if}
                    \State{$ S_* \gets S $} \label{alg:grasp:updateSs}
                    \State{$ f_* \gets f(S_*) $} \label{alg:grasp:updatefs}
                \EndIf
            \EndWhile \label{alg:grasp:forend}
            \State \Return $S_*$ \label{alg:grasp:return}
        \end{algorithmic}
\end{algorithm}

\subsection{Path Relinking}
\label{sec:PR}

Path relinking (PR) was motivated by a concept called proximate optimality principle, introduced in Chapter 5 of \citet{glolag1997}. The proximate optimality principle (POP) may be characterized as a heuristic counterpart of the so called Principle of Optimality in dynamic programming.  However, POP does not incur in the curse of dimensionality produced by the explosion of state variables that  occurs when dynamic programming is applied to combinatorial problems. The proximate optimality principle stipulates that good solutions at one level are likely to be found close to good solutions at an adjacent level. In GRASP and path relinking, ``level" refers to a stage in the constructive process. The intuition is to seek improvements at a given level before moving to the next level.

The POP notion implies a form of connectivity for the search space that path relinking seeks to exploit. PR trajectories, which are guided by elite solutions, are likely to go through regions where new elite solutions reside, as long as appropriate neighborhood definitions are used.  The result is a type of focused diversification that complements the “sampling” nature of GRASP. In this context, an elite set consists of solutions, found during the search, that are the best according to their objective function value. In the original PR proposal, upon identifying a collection of one or more elite solutions to guide the path of an initiating solution, the attributes of these guiding solutions were assigned preemptive weights as inducements to be selected. Larger weights were assigned to attributes that occur more often in the set of guiding solutions, allowing bias to give increased emphasis to solutions with higher quality or with special features. However, most PR implementations select one solution from the elite set to act as the sole guiding solution.

The elite set in the original hybridization of GRASP with PR \citep{lag1999} consists of three solutions. Therefore, the procedure starts with three GRASP iterations to populate the elite set. From the fourth iteration until the last, the  local optimum at the end of the improvement phase is used as the initiating solution for a PR phase. The guiding solution is randomly chosen from the elite set. At the end of the PR phase, the elite set is updated if the best solution found during the path relinking has a better objective function value than the worst elite solution. The PR includes a neighborhood search that is applied after every few steps of the process. The number of steps between each application of the neighborhood search is controlled by a search parameter that attempts to find an effective separation between every pair of neighborhood searches while taking into consideration the additional computational effort. 

In Section \ref{sec:AdvancedPR}, we describe alternative PR designs, such as greedy, greedy randomized, and mixed PR. Section \ref{sec:GRASP_PR} revised GRASP with  PR designs, including dynamic, static and evolutionary variants. One of the main contributions of the paper comes in Section \ref{sec:history}, where a historical review and critical analysis of the variants described in the previous sections is presented. We divide the evolution of the methodology into two decades, the creative decade from 2000 to 2010, and the consolidation decade, from 2010 to the present. Section \ref{sec:other_methods} describes the main connections of GRASP with PR with other metaheuristics. The paper finishes with Section \ref{sec:comput} with an extensive computational experimentation with two well-known hard optimization problems: the Linear Ordering Problem and the Max-Cut problem, to help us to evaluate the performance of the different variants and designs reviewed. The associated conclusions end our journey on these 20 years of GRASP with PR.

\section{Advanced path relinking strategies}\label{sec:AdvancedPR}

\begin{figure}[ht]
	\centering
 \includegraphics[width=0.35\textwidth]{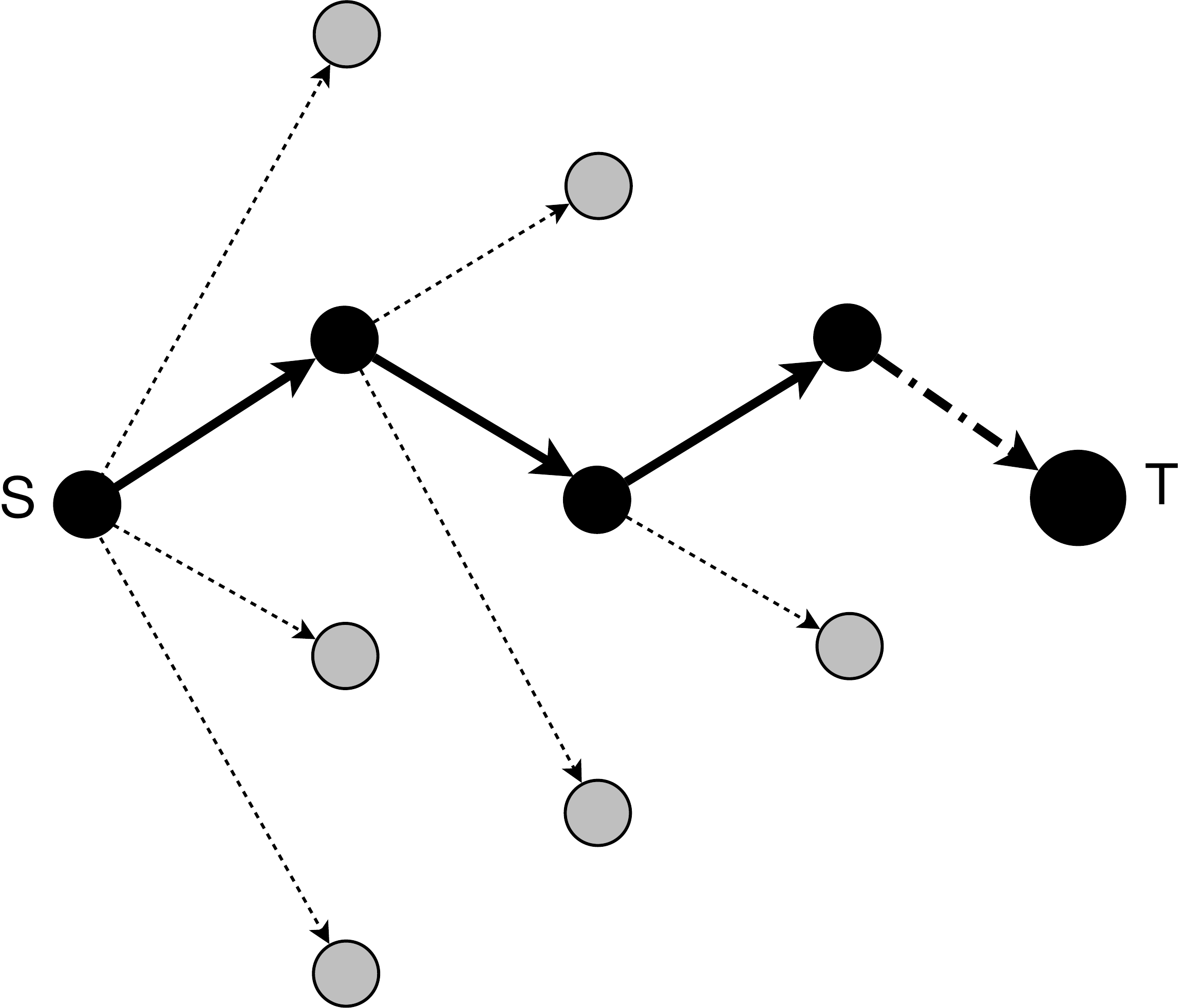}
	\caption{Forward path-relinking: a path is traversed
from the initial solution $S$ to a guiding solution $T$ at least as good as $S$.}
		\label{fig:fpr}
\end{figure}

In Section~\ref{sec:origins} we considered the original variant of path relinking proposed in \cite{lag1999}. This variant, which we refer to here as \textit{greedy forward path relinking}, or simply forward path relinking, uses as the initiating (or initial) solution $S$ the output of the local search phase of GRASP and as guiding solution $T$ one selected from the elite set.  More generally, we can say that given two solutions $S$ and $T$ to relink (with $f(T) \ge f(S)$), forward path relinking uses the better solution, i.e. $T$ for a maximization problem, as the guiding solution and other ($S$) as the initiating solution. This variant is referred to as greedy since upon incorporating in $S$ the attributes of $T$ not present in $S$ the updated initiating solution $S'$ is selected as the one resulting from the introduction of the attribute that leads to the solution with best cost. Figure \ref{fig:fpr} illustrates the path from $S$ to $T$ with its \textit{intermediate solutions}, depicted in black, and several solutions in their neighborhoods, depicted in gray. 

Path relinking can be thought of as a constrained neighborhood search, where the search is limited to explore the solutions in the neighborhood with characteristics of the guiding solution. The selected neighborhood will determine the set of solutions visited by path relinking.  For example, consider a solution represented by a permutation and consider two neighborhoods, swap and insert.  In swap $(\pi^1, \ldots, \pi^i, \ldots, \pi^j, \ldots, \pi^n)$ and $(\pi^1, \ldots, \pi^j, \ldots, \pi^i, \ldots, \pi^n)$ are neighbors because $\pi^i$ and $\pi^j$ swap their positions whereas in insert $(\pi^1, \ldots, \pi^{i-1}, \pi^{i}, \ldots, \pi^n)$ and $(\pi^1, \ldots, \pi^{i-1}, \pi^j, \pi^{i}, \ldots, \pi^n)$ are neighbors since $\pi^j$ is inserted in position $i$.

\begin{figure}[ht]
		\centering
		\includegraphics[width=0.5\textwidth]{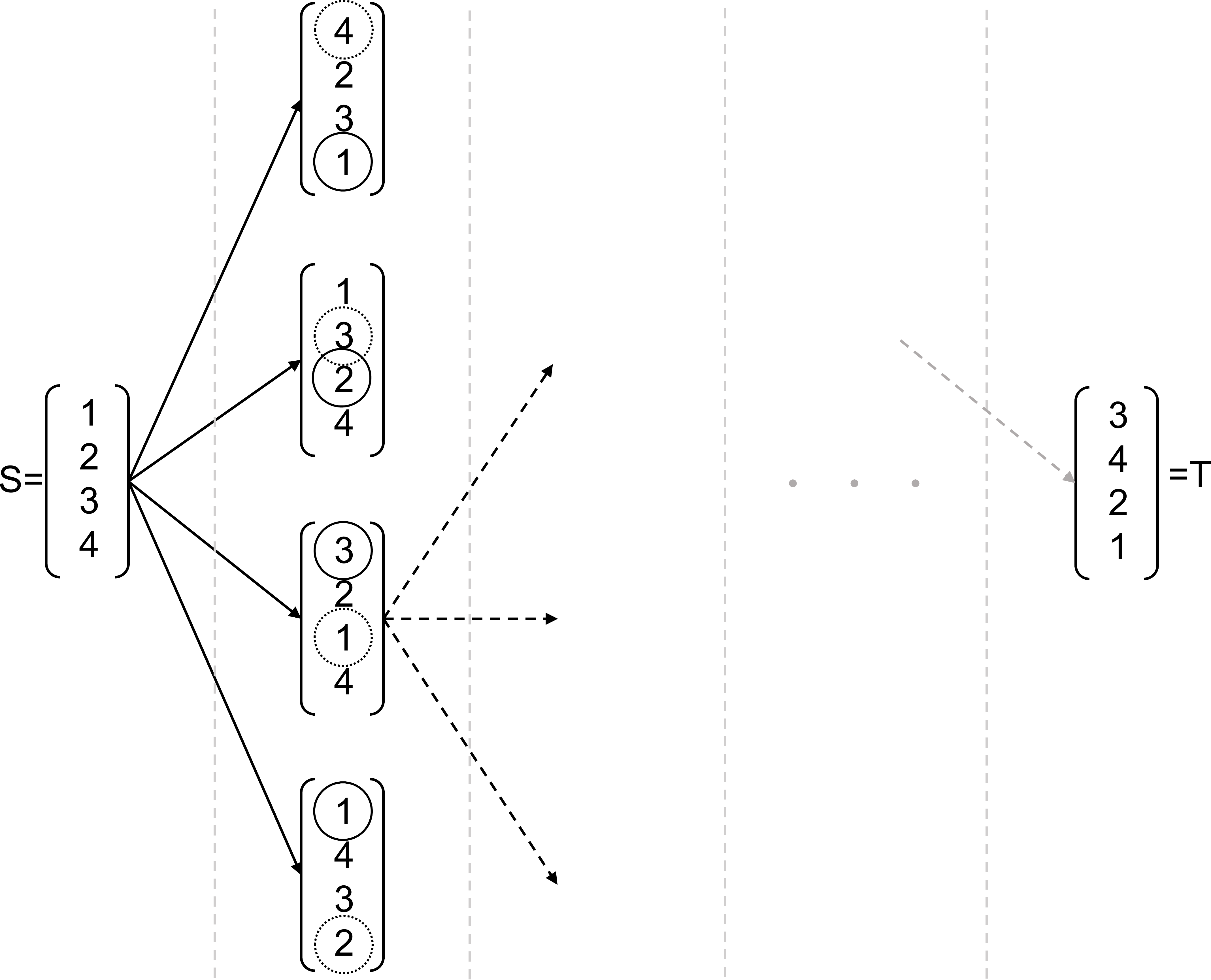}
	\caption{Intermediate solutions in a swap neighborhood.}
	\label{fig:neig.swaps}
\end{figure}

\begin{figure}[ht]
		\centering
    \includegraphics[width=0.5\textwidth]{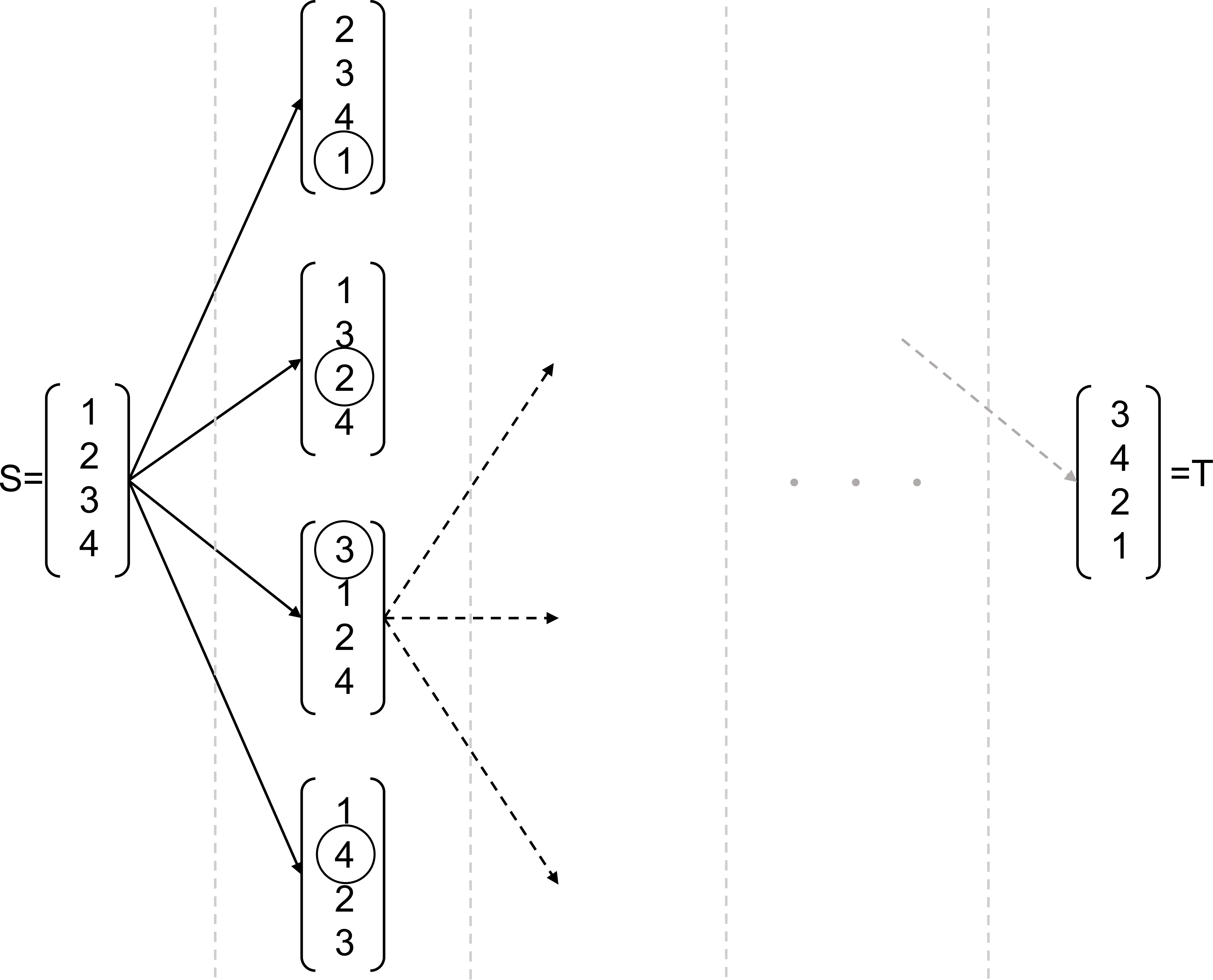}
	\caption{Intermediate solutions in an insertion neighborhood.}
	\label{fig:neig.insert}
\end{figure}

Figures \ref{fig:neig.swaps} and \ref{fig:neig.insert} illustrate the application of PR in an example in which solutions are represented with $4$-element permutations. In particular, Figure \ref{fig:neig.swaps} shows the neighborhood of the initiating solution $S=\{1,2,3,4\}$ when applying swap moves to reach the target solution $T=\{3,4,2,1\}$. This neighborhood contains four solutions obtained by exchanging each element with the corresponding element in the position given by the target solution. For example, element $1$ is in position $1$ in the initial solution $S$, and it is in position $4$ in the target solution $T$. Then, PR generates a solution in which element $1$ is in its position in $T$. Note however that since we are considering a swap move, this implies to swap element $1$ with the element that is in position 4 in $S$, which is element 4 in this case. Therefore, this solution effectively swaps elements 1 and 4 in $S$, obtaining the intermediate solution $\{4,2,3,1\}$. To highlight that the origin of this swap move is to place element 1, it is depicted in a solid circle in the figure. On the other hand, element 4 is depicted with a dashed circle to represent that it was also moved to perform the swap of 1. We obtain in a similar way the other three intermediate solutions depicted in Figure \ref{fig:neig.swaps}. 

Figure \ref{fig:neig.insert} represents a PR implementation based on an insertion move that is simpler than the swap move. In this case, each element in the initiating solution $S$ is directly inserted in its position in the target solution $T$, and the rest of the elements are simply shifted. Comparing Figures \ref{fig:neig.swaps} and \ref{fig:neig.insert} is interesting to observe the different intermediate solutions that they generate in the path from $S$ to $T$.

\begin{figure}[ht]
	\centering
 \includegraphics[width=0.35\textwidth]{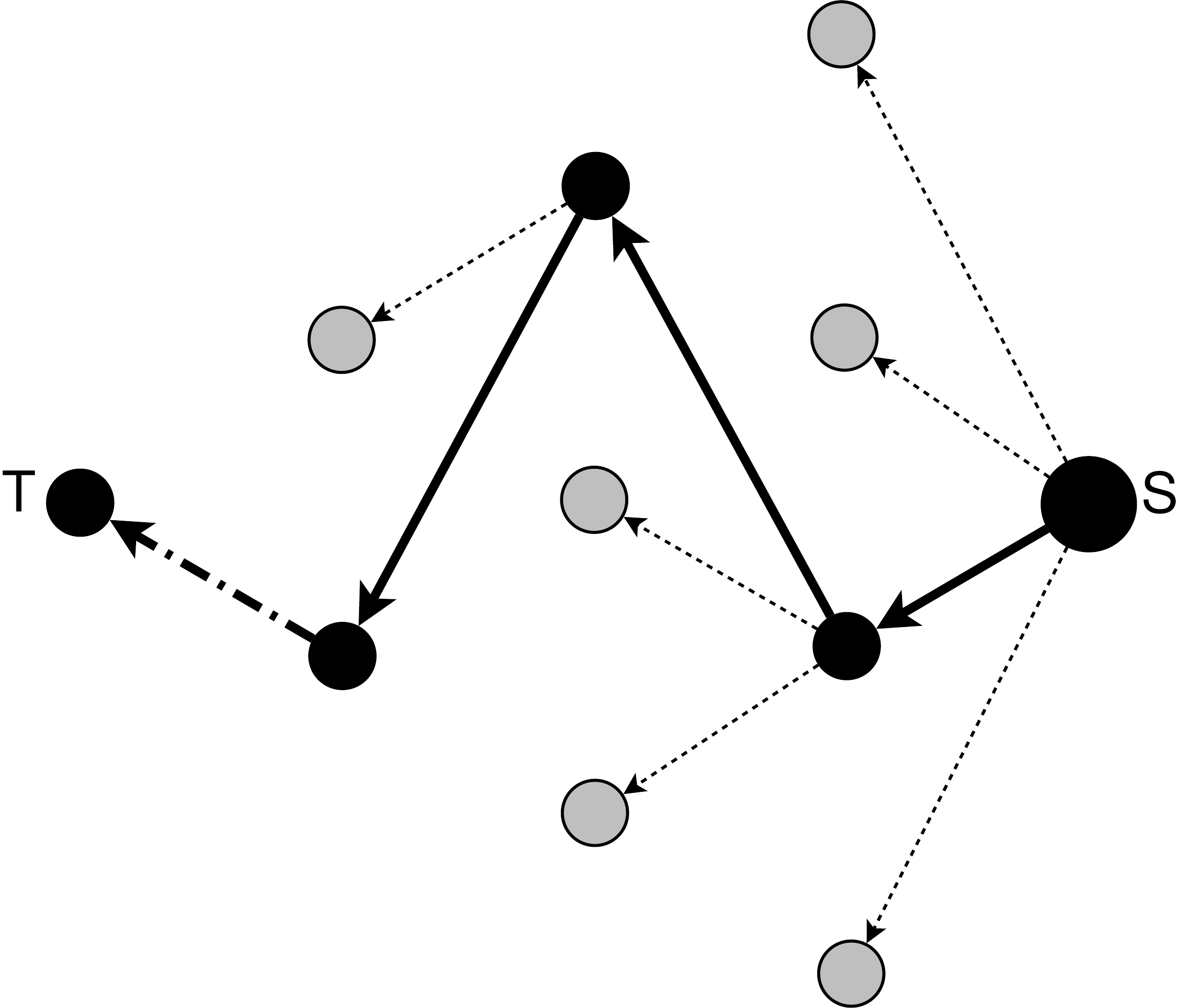}
	\caption{Backward path-relinking: a path is traversed from the initial solution $S$ to a guiding solution $T$ that is not better than $S$.}
		\label{fig:bpr}
\end{figure}

The \textit{backward path relinking} between $S$ and $T$ is similar to its forward variant except that now the initiating solution is $S$, where $f(S) \ge f(T)$ and the guiding solution is $T$. Figure \ref{fig:bpr} illustrates it. In \textit{back and forward path relinking} one first applies backward path relinking followed by forward path relinking.

\subsection{Mixed path relinking}

In forward and in backward path relinking only one of the constrained neighborhoods of $S$ or $T$ is fully explored, as is illustrated in Figures \ref{fig:fpr} and \ref{fig:bpr}, respectively. For example, in forward path relinking only the constrained neighborhood of the initiating solution $S$ is fully explored.  That is because in each round of path relinking one element of the guiding solution is fixed in each intermediate solution and as the rounds proceed the size of the constrained neighborhoods decreases by one element per round. 

Consider as an example the permutation solution of four elements where the initiating solution is $S = \{ 1, 2, 3, 4\}$ and the guiding solution is $T = \{ 3, 4, 2, 1 \}$.  In the first round, the constrained neighborhood has four solutions corresponding to the insertion of each element of $S$ in its corresponding position in $T$.  Suppose the solution with best cost (greedy solution) corresponds to placing element 3 of $S$ in position 1 as imposed by the guiding solution. This results in $S_1 = \{ \underline{3}, 1, 2, 4 \}$.  From this point on element 3 is fixed to position 1 in the intermediate permutations, and at this point the size of the constrained neighborhood is reduced from 4 to 3 since we can now only move elements 1, 2, or 4 to their positions in the guiding solution. This way if we next set element 4, then $S_2 = \{ \underline{3}, \underline{4}, 1, 2 \}$ with constrained neighborhood of size 2. If the next element selected is element 2, we reach $S_3 = \{ \underline{3}, \underline{4}, \underline{2}, 1 \}$ whose constrained neighborhood only contains the guiding solution.

\begin{figure}[ht]
	\centering
 \includegraphics[width=0.4\textwidth]{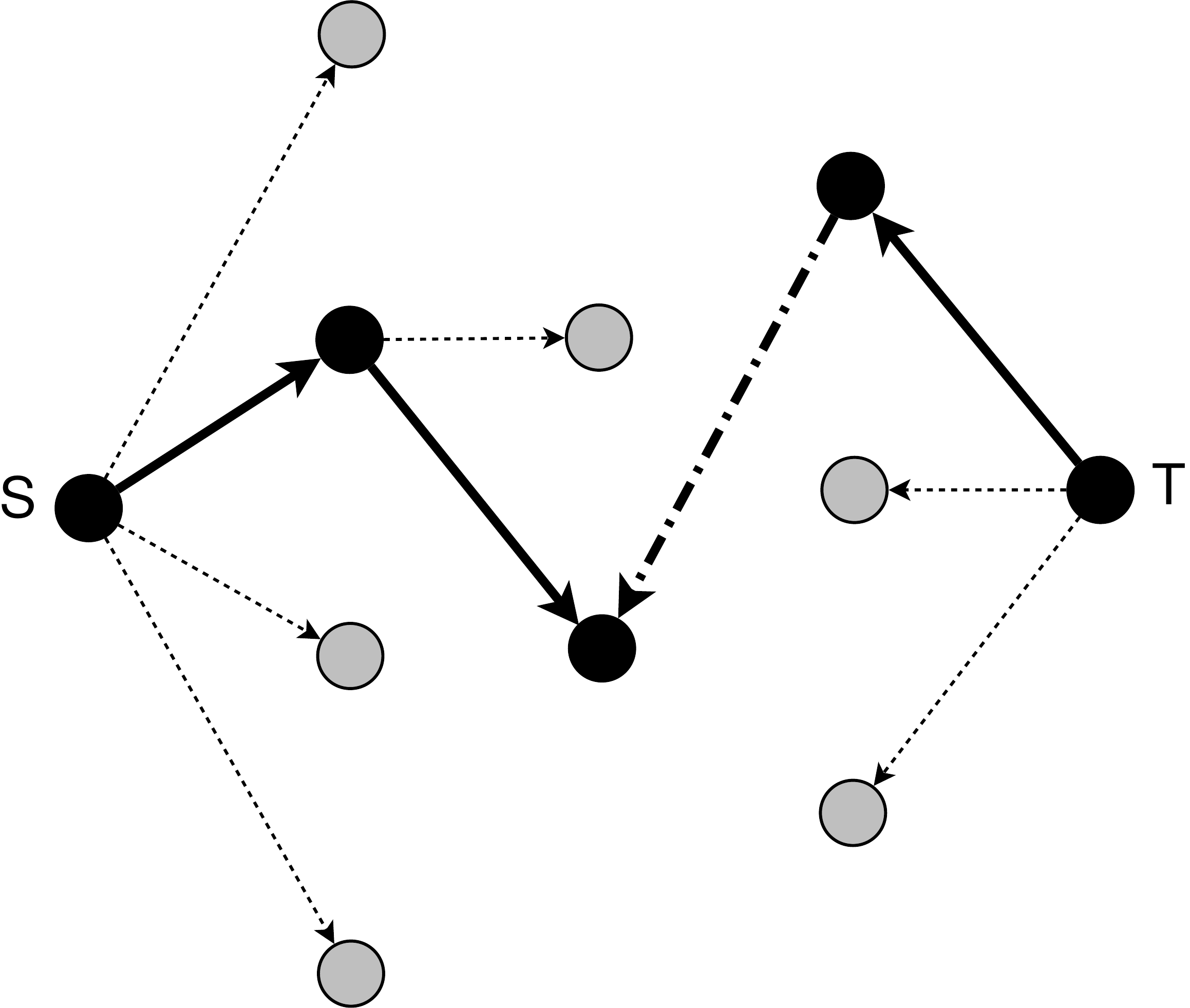}
	\caption{Mixed path-relinking: two subpaths are traversed, one starting at $S$ and the other at $T$,
which eventually meet in the middle of the trajectory connecting $S$ and $T$.}
		\label{fig:mpr}
\end{figure}

One way to explore the neighborhoods of both $S$ and $T$ is to apply \textit{mixed path relinking}, as illustrated in Figure \ref{fig:mpr}. In this variant forward path relinking is applied but at each round the roles of the guiding and the initiating solutions are reversed.  In this way the algorithm explores the constrained neighborhoods of both the initial guiding and initiating solutions. 

Consider again the example where the initiating solution is $S = \{ 1, 2, 3, 4\}$ and the guiding solution is $T = \{ 3, 4, 2, 1 \}$. The first round leads us to the same intermediate solution as in the above example of forward path relinking, i.e.  $S_1 = \{ \underline{3}, 1, 2, 4 \}$.  In the next round we exchange the roles of $S_1$ and $T$, i.e. $T \gets S_1$ and $S \gets T$ and again apply path relinking from $S$ to $T$.
At this point we have three choices, i.e. elements 1, 2, or 4. This way mixed path relinking explores the neighborhood of the initial guiding solution.
If we move element 4 in $S = \{ 3, 4, 2, 1 \} $ to its position in $T = \{ \underline{3}, 1, 2, 4 \}$, this will result in $S_2 = \{ \underline{3}, 2, 1, \underline{4} \}$.

\subsection{Path cardinality and truncated path relinking} \label{sec:truncated}

At each round of path relinking an element (or attribute) of the guiding solution is added to the initiating solution making the updated solution one step closer to the guiding solution.  Let $\Delta$ be the set of indices of the attributes present in the guiding solution $T$ that are not present in the initiating solution $S$.
The number of elements to be added is therefore the cardinality of the symmetric difference $|\Delta|$. For example, if a solution is represented by a 0-1 indicator vector of size $n$, then $\Delta = \{ j = 1,2, \ldots, n : S_j \not= T_j\}$.  If a solution $S$ is represented by a permutation $\pi^S = (\pi^S_1, \pi^S_2, \ldots, \pi^S_n)$ of $n$ integers, then $\Delta = \{ j = 1,2, \ldots, n : \pi^S_j \not= \pi^T_j\}$.
The path cardinality is therefore $|\Delta|-1$ since at the last round the solution visited is a neighbor of the guiding solution $T$.

\begin{figure}[ht]
	\centering
 \includegraphics[width=0.6\textwidth]{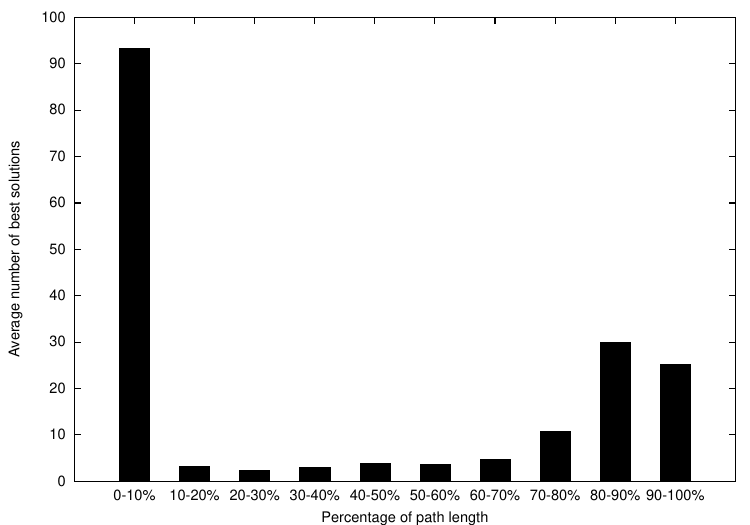}
	\caption{Histogram of the average of the best solution found in the path \citep{res2010}.}
		\label{fig:trunc}
\end{figure}

It is a reasonable assumption that the most promising neighborhoods to search for improving solutions are the neighborhoods of the initiating solution $N(S)$ and of the guiding solution $N(T)$.  The neighborhoods of solutions within these neighborhoods may also be promising, but as one moves further from them in the path they become less promising.
A strategy to take advantage of this observation is \textit{truncated path relinking} where only solutions on the path close to $S$ or $T$ are visited.  In forward path relinking only solutions near the initiating solution are visited while in backward path relinking only solutions near the guiding solution are visited. In back and forward as well as in mixed path relinking solutions near both the initiating and guiding solutions are visited.
As an example of above argument, consider the GRASP with path relinking for max-min diversity described in
\citet{res2010} a plot (shown here in Figure~\ref{fig:trunc}) shows that the largest number of best solutions found by path relinking are indeed found near $S$ and $T$. The figure shows the average number of best solutions found in the first $10\%$ of the path, the number of best solutions in the second $10\%$ of the path (from 10\% to 20\%), and so on. The figure confirms the hypothesis that the best solutions are mainly obtained at the beginning of the path. However, good solutions are also obtained in the final part of the path.

Truncated path relinking is a good mechanism to speed up the algorithm. However, the larger the cardinality path, the better the performance of the truncated variant relative to the nontruncated variant. 

\subsection{Minimum distance required for path-relinking}
\label{s_distance}


It is natural that two solutions that are to be relinked be locally optimal.
Let us assume that we want to relink two locally optimal
solutions $S$ and $T$.
If these two solutions differ 
by only one of their components,
then the path directly connects the two solutions and no solution,
other than $S$ and $T$, is visited.

Since $S$ and $T$ are both local optimal, then $f(S) \geq f(R)$ for all $R \in N(S)$ and
$f(T) \geq f(R)$ for all $R \in N(T)$.
If $S$ and $T$ differ by exactly two elements, then
any path between $S$ and $T$ visits exactly
one intermediate solution $R \in N(S) \cap N(T)$.
Consequently, solution $R$ cannot be better than either $S$ or $T$.

\begin{figure}[ht]
\centering
\includegraphics[width=0.6\textwidth]{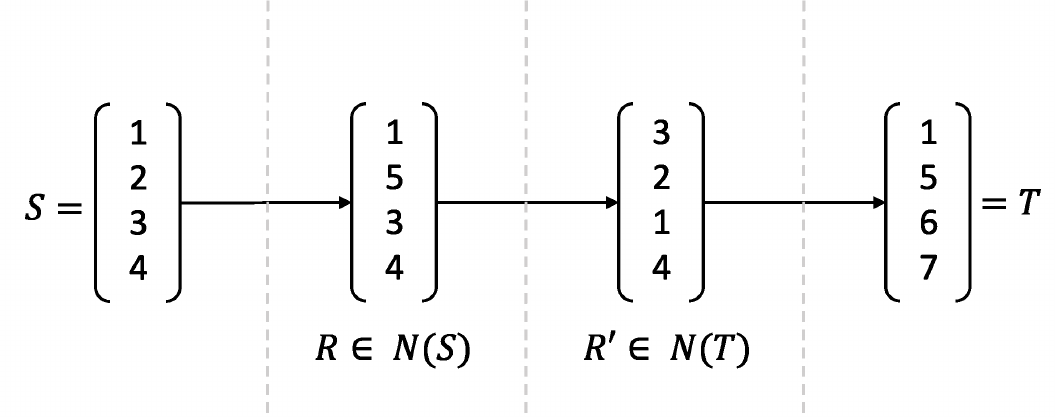}
\caption{Path between two solutions that differ in three elements. 
}
\label{fig:dist-3}
\end{figure}

Likewise, if
$S$ and $T$ differ by exactly three elements, then
any path between them visits two
intermediate solutions $R \in N(S)$ and
$R' \in N(T)$ and, consequently, neither $R$ nor $R'$ can be
better than both $S$ and $T$. Figure~\ref{fig:dist-3} illustrates this  case in an example with $S=\{1,2,3,4\} $ and $T=\{1,5,6,7\}$, where the solutions are represented with the four selected elements out of eight in a selection problem. These two local optima, $S$ and $T$, differ by three elements, and therefore the path between them contains at most two intermediate solutions, say $R$ and $R'$. The first intermediate solution, $R=\{1,5,3,4\}$, is in the neighborhood of $S$, and considering that $S$ is a local optimum, $R$ cannot be better than $S$. Similarly,  $R'$ cannot be better than $T$. Therefore, it is only possible to find an improving
solution in the path between $S$ and $T$ 
when these two solutions differ by at least four elements. In other words, the minimum distance to apply path relinking can be set when the two relinked solutions differ in at least four elements (in a selection problem). Generalizing this argument, we expect to find better solutions when applying path relinking to significantly different solutions that generate longer paths.

\subsection{Multiple parents in Path Relinking}

\begin{figure}[ht]
\centering
\includegraphics[width=0.4\textwidth]{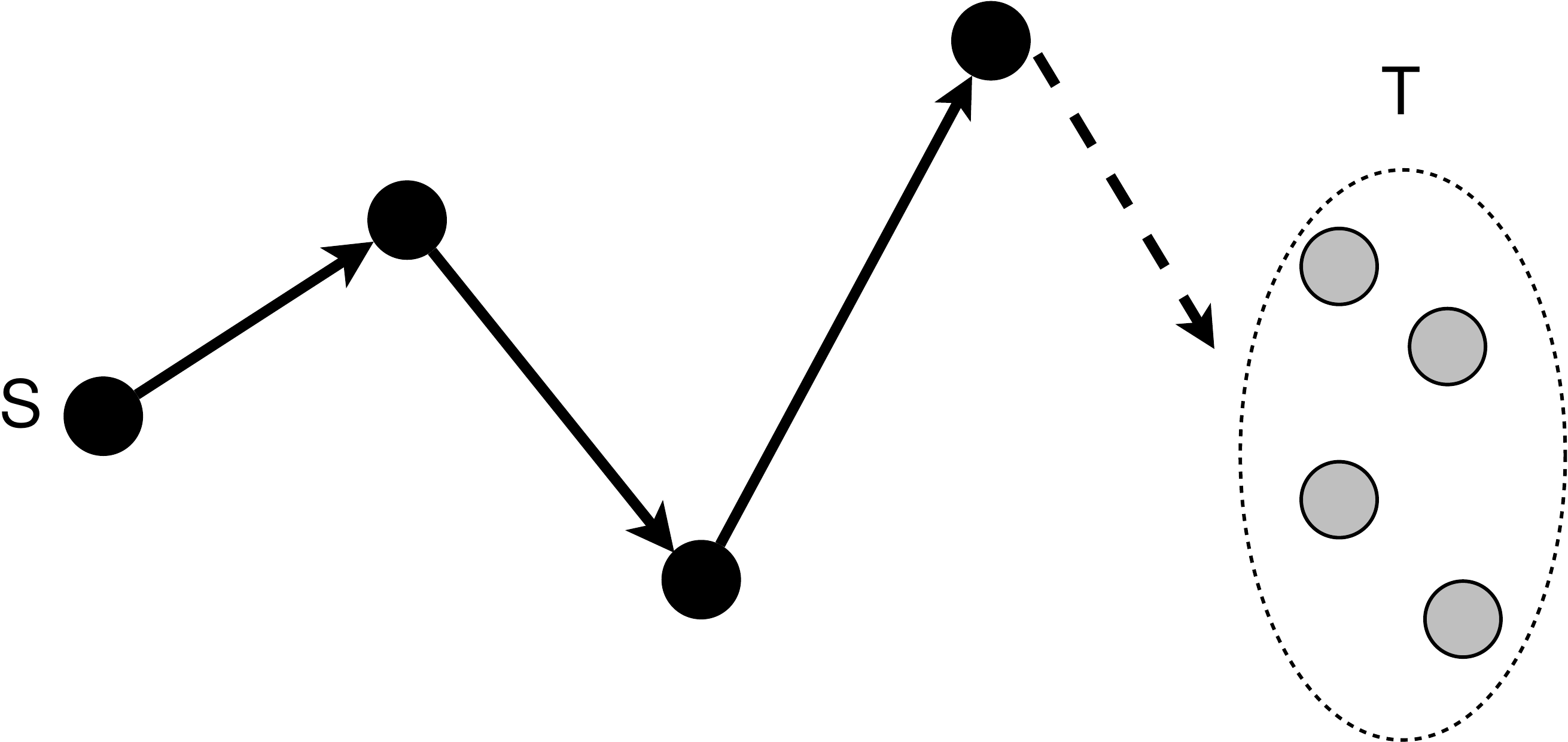}
\caption{Multiple guiding solutions. 
}
\label{fig:multiplePR}
\end{figure}

In previous variants of path relinking, a path is generated by going from an initial solution $S$ to a guiding solution $T$. Instead of that, the path can be created considering a collection of guiding solutions $\mathcal{T}$, as illustrated in Figure \ref{fig:multiplePR}. From an intermediate solution $S_i$, the next solution $S_{i+1}$ is obtained by incorporating an attribute from the union of attributes of the guiding solutions in  $\mathcal{T}$. This design can be interpreted as approaching, from the initiating solution, to the region where the guiding solutions belong. Note that the path does not necessarily match one of the guiding solutions, but it can terminate when a certain criterion is satisfied (based on the number of attributes incorporated to the intermediate solutions).

In advanced designs, each attribute is counted (weighted) in accordance to the number of times that it appears in the guiding set of solutions. This strategy was first proposed in \cite{glo1994} in the context of the Tabu Search methodology. The idea behind the use of multiple guiding solutions is to preserve desirable solution properties, such as feasibility in some optimization problems, without requiring the application of a repair procedure. Due to its complexity has been rarely applied in combination with GRASP.

\subsection{Interior and Exterior Path Relinking}

Most PR implementations only consider the between-form of PR (Interior Path Relinking). This section discusses the beyond-form of path relinking, introduced in \citet{Glover2014} and called \textit{Exterior Path Relinking} (EPR). 
Instead of introducing into the initiating solution characteristics present in the guiding solution, EPR introduces in the initiating solution 
characteristics not present in the guiding solution.
Specifically, it removes from the initiating solution those
elements or attributes which also belong to the guiding solution, obtaining
intermediate solutions which are further away from both the
initiating the guiding solutions.

\begin{figure}[ht]
\centering
\includegraphics[width=0.6\textwidth]{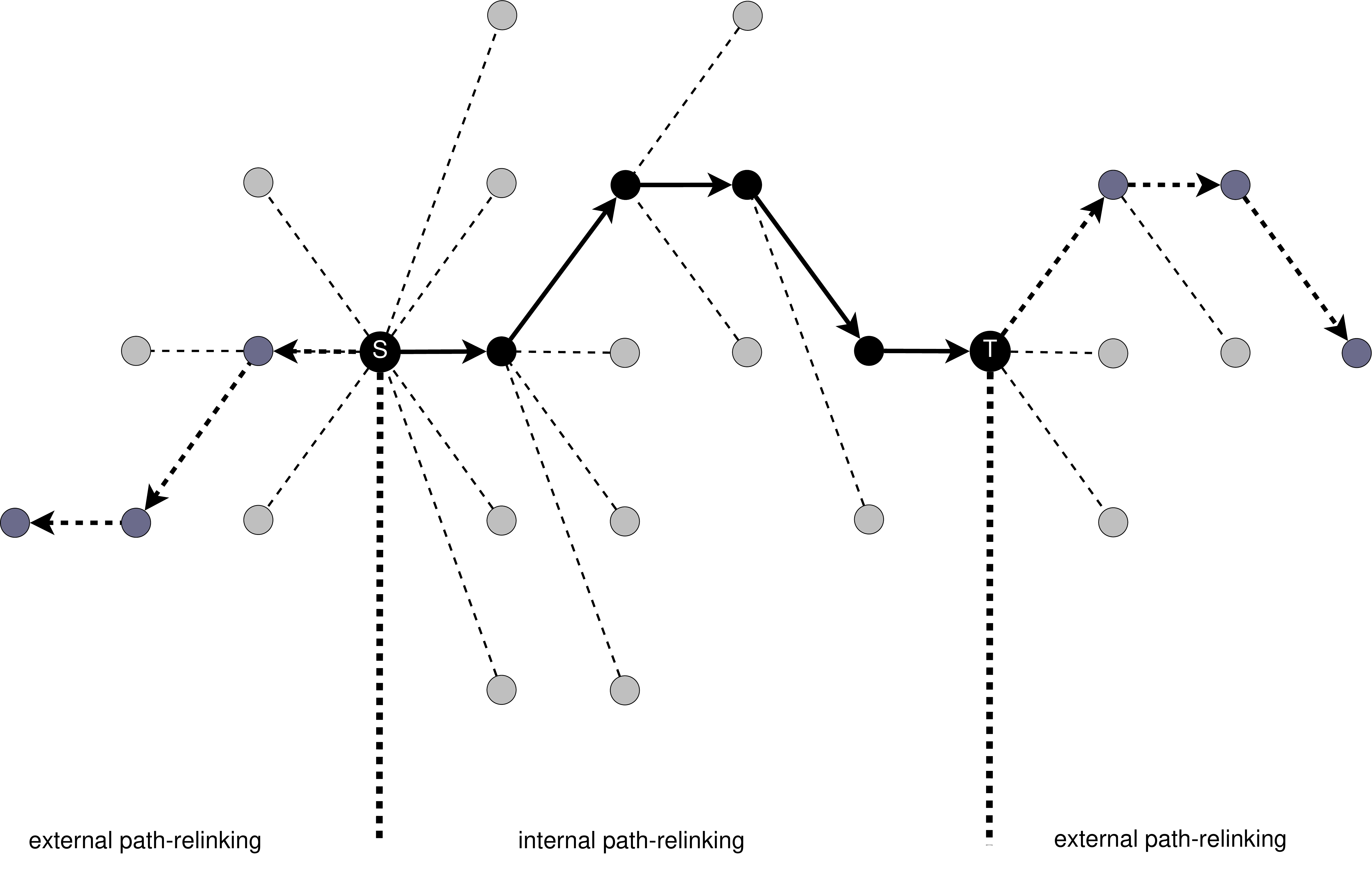}
\caption{Internal and external path-relinking.
}
\label{fig:epr}
\end{figure}

Figure \ref{fig:epr}, adapted from \citet{res2016}, illustrates the application of both interior (also called internal) and exterior (external) PR. An internal path (black arcs, black nodes) from solution $S$ to solution $T$ and two external (black hashed arcs, dark gray nodes) paths, one emanating from solution $S$ and the other from solution $T$.

The relevance of paths that go beyond the initiating and guiding solutions was broached in \citet{glo1997} as follows: The scope of strategies made available by path relinking is significantly affected by the fact that the term neighborhood has a broader meaning in tabu search than it typically receives in the popular literature on search methods.  Often, the neighborhood terminology refers solely to methods that progressively transform one solution into another.  Such neighborhoods are called transition neighborhoods in tabu search, and are considered as merely one component of a collection of neighborhoods that also include those operating in regions beyond solutions previously visited.

\subsection{Greedy randomized adaptive path relinking }\label{grpr}

As aforementioned, path relinking traverses a path in the solution space linking the initiating solution to the guiding solution. However, the number of paths linking these two solutions grows exponentially with the cardinality of the symmetric difference between the two extreme solutions. Standard path relinking approaches only search one of these paths.  The path that is selected is what we call the 
\textit{greedy path}.  This is so because at each step in path relinking all solutions obtained by inserting into the current solution the different attributes of the guiding solution not found in the current solution are evaluated and the solution with the best cost (the greedy choice) is selected as the new current solution. 

\textit{Greedy randomized adaptive path relinking (GRAPR)} \citep{far2005} enables the exploration of many good paths linking the two extreme solutions. As with the standard path relinking, GRAPR evaluates the consequence of the insertion of each attribute of the guiding solution not present in the current solution.  However,  instead of making the greedy choice as is done in standard path relinking, GRAPR builds a restricted candidate list (RCL) of near-greedy insertions and then selects a candidate from the RCL at random. This process results in a \textit{semi-greedy path}. As with GRASP, the amount of randomness or greediness of this process can be controlled by fixing the cardinality of a cardinality-based RCL or changing the acceptance cutoff in a value-based RCL.  While it is possible that GRAPR and a standard path relinking explore the same path, it become more and more unlikely as the RCL grows in size.

GRAPR can be implemented as different variants as is the case for PR, such as forward, backward, back and forward, mixed, and/or truncated.
It can be incorporated into a GRASP as a multi-start procedure in place of each application of PR done in standard GRASP with PR.  Furthermore, it can be incorporated into \textit{evolutionary PR} as we will see in Section~\ref{s_evPR}.

\section{GRASP with Path Relinking designs}
\label{sec:GRASP_PR}

In this section we review different alternatives and strategies in which GRASP and Path Relinking have been hybridized and implemented. Path Relinking makes use of an elite set (ES) of good and diverse solutions that are collected over the iterations of the GRASP with PR algorithm. The use of an elite set adds a memory mechanism to GRASP, since the best solutions across iterations are kept in the ES. From a methodological point of view, this makes GRASP with PR a hybrid method that brings together memory-based and memory-less methods.

\subsection{Dynamic PR}\label{dynamic}

Two strategies to manage the pool of solutions in ES are \textit{dynamic} PR and \textit{static} PR. In the former, as shown in line~\ref{alg:alg-2:pr} of Algorithm~\ref{alg:alg-2}, path relinking is applied after each local search phase of GRASP. In this type of PR, an element $S'$ of the ES is relinked with the locally optimal solution $S$ found during the local search phase of GRASP. In the latter, as shown in lines~\ref{alg:alg-3:for} to~\ref{alg:alg-3:endfor} of Algorithm~\ref{alg:alg-3}, path relinking is applied to all pairs of elite set solutions after the conclusion of the GRASP iterations.

This \textit{dynamic} strategy to update the elite set was introduced in \citet{lag1999}. Solution 
$S'$ in \textrm{ES} is selected uniformly at random. In some implementations, this choice is made at random with probability proportional to the cost \citep{res2010} or to the symmetric difference  \citep{res2004}. In this design, \textrm{ES} is first populated with the first $k$ GRASP solutions. In the original proposal, $k$ was set to 3. In line~\ref{alg:alg-2:kisone} of Algorithm~\ref{alg:alg-2}, $k$ is set to 1, but recent studies have applied larger values (between 10 and 20). The local search method is usually applied to the output of PR, and the resulting solution is directly tested for inclusion in \textrm{ES}, as shown in line~\ref{alg:alg-2:updateES}. If successful, it can be used as guiding solution in later applications of PR.

\begin{algorithm}[htbp]
	\caption{\texttt{Dynamic GRASP with PR}}
	\label{alg:alg-2}
	\begin{algorithmic}[1]
            \State{$\mathit{PRflag} \gets \mathtt{false}$} \label{alg:alg-2:es}
           \State{$f_* \gets -\infty$} \label{alg:alg-2:fs}
            \While {\textrm{stopping criterion is not satisfied}} \label{alg:alg-2:while}
                \State{$ S \gets \texttt{Greedy\_Randomized\_Adaptive\_Construction}(\alpha) $} \label{alg:alg-2:grc}
                \State{$ S \gets \texttt{Local\_Search}(S) $} \label{alg:alg-2:ls}
                \If{$\mathit{PRflag}$}
                    \State{$S'\gets \texttt{Select\_Solution}(ES)$} \label{alg:alg-2:select}
                    \State{$S \gets \texttt{PR}(S', S)$} \label{alg:alg-2:pr}
                    \State{$  ES \gets \texttt{Update\_Elite}(ES,S) $}
                \Else
                    \State{$\mathit{PRflag} \gets \mathtt{true}$}
                    \State{$\mathit{ES} \gets \{ S \}$} \label{alg:alg-2:kisone}
                \EndIf
                \
                \If{$f(S) > f_*$} \label{alg:alg-2:updateES}
                    \State{$ S_* \gets S $} \label{alg:alg-2:updateSs}
                    \State{$ f_* \gets f(S_*) $} \label{alg:alg-2:updatefs}
                \EndIf
            \EndWhile \label{alg:alg-2:endWhile}
            \State \Return $S_*$ \label{alg:alg-2:return}
        \end{algorithmic}
\end{algorithm}

\subsection{Static PR}\label{static}

In the \textit{static} design, PR is applied as a post-processing after the application of GRASP. As shown in Algorithm~\ref{alg:alg-3}, the method is called \textit{static} because it consists of two separate phases, a GRASP phase and a PR phase where the elite set does not change. In this way, PR does not start until the best GRASP solutions have been identified in the first phase, and therefore only the best solutions are used to guide the path in PR. As in the \textit{dynamic} strategy, the local search method is applied to the output of PR, but here this newly generated solution is not submitted to ES for potential inclusion. The PR phase terminates when all pairs of solutions in ES have been relinked with PR.

\subsection{Elite set management}\label{elite}

The elite set $\mathrm{ES}$ used in GRASP with PR consists of a fixed number of diverse high-quality solutions found during the iterations of the algorithm.  Elite set management consists in the initialization
of the elite set, selection of a solution from the elite set, and updating the elite set every time a new improving solution is found.

When the initial elite set is constructed, a newly found solution is added to it if this solution is sufficiently different from all solution currently in the elite set.  If the elite set is still empty, then the solution on hand is added to the elite set as its first solution.  

Once the elite set is fully populated and a new solution $S$ with cost $c(S)$ is found, this solution will be added to the elite set if it is better than the worst solution currently in the elite set.  Since the size of the elite is fixed, if a new solution is added to the set, then necessarily an elite set solution will need to be removed.  The goal is to improve the average solution value of the elite set while increasing its diversity or decreasing it as as little as possible.
To do this, \citet{res2004} suggest selecting for removal the solution $T^*$ most similar to $S$ among those elite set solutions of worse quality, i.e. 
\begin{displaymath}
T^* = \mathrm{argmin} \{ \Delta(T,S) \;\; \forall T \in \mathit{ES} \;|\; c(T) < c(S) \}.
\end{displaymath}

During each iteration of the dynamic GRASP with PR, a solution $T$ is chosen from the pool and relinked with $S$, the most recently identified GRASP solution. 
Early implementations of GRASP with PR
proposed the random selection of a solution $T$ with equal probability \citep{lag1999}. However, this can often result in picking a solution too akin to $S$, reducing the chances of discovering improved new solutions. To address this issue, \citet{res2004} select a solution $T$ from the pool with probability proportional to the symmetric difference $\Delta(T,S)$ with respect to $S$. 

\begin{algorithm}[htbp]
	\caption{\texttt{Static GRASP($\alpha$) with PR}}
	\label{alg:alg-3}
	\begin{algorithmic}[1]
            \State{$ES \gets \emptyset$} \label{alg:alg-3:es}
           \State{$f_* \gets -\infty$} \label{alg:alg-3:fs}
           
            \While {\textrm{stopping criterion is not satisfied}} \label{alg:alg-3:while}
                \State{$ S \gets \texttt{Greedy\_Randomized\_Adaptive\_Construction}(\alpha) $} \label{alg:alg-3:grc}
                \State{$ S \gets \texttt{Local\_Search}(S) $} \label{alg:alg-3:ls}
                \State{$  ES \gets \texttt{Update\_Elite}(ES) $} \label{alg:alg-3:updateES}
                \If{$f(S) > f_*$} \label{alg:alg-3:if}
                    \State{$ S_* \gets S $} \label{alg:alg-3:updateSs}
                    \State{$ f_* \gets f(S_*) $} \label{alg:alg-3:updatefs}
                \EndIf
            \EndWhile \label{alg:alg-3:endwhile}
            \For{$ S, T \in ES$} \label{alg:alg-3:for}
                \State{$S \gets \texttt{PR}(S, T)$} \label{alg:alg-3:pr}
                \If{$f(S) > f_*$} \label{alg:alg-3:ifPR}
                    \State{$ S_* \gets S $} \label{alg:alg-3:updateSs2}
                    \State{$ f_* \gets f(S_*) $} \label{alg:alg-3:updatefs2}
                \EndIf
            \EndFor \label{alg:alg-3:endfor}
            \State \Return $S_*$ \label{alg:alg-3:return}
        \end{algorithmic}
\end{algorithm}

\subsection{Evolutionary Path Relinking}
\label{s_evPR}

 \citet{res2004} proposed Evolutionary Path Relinking, where the solutions in the elite set are evolved in a similar way that the reference set evolves in scatter search \citep{lag2003}.
The GRASP with PR iterations yield multiple local optima that are typically close in quality to the best locally optimum solution found. 
Evolutionary path relinking, often implemented as a post-processing phase, combines these solutions with the goal of deriving even better solutions. This procedure uses the elite solution set, elaborated during the dynamic GRASP with PR iterations, as its input. Each pair of solutions from this set is combined with through path-relinking. The resulting solutions are collated into a fresh set of elite solutions, adhering to the guidelines given in Section~\ref{elite}, resulting in a new iteration of the evolutionary algorithm. The algorithm continues until a cycle is reached where there is no improvement from the prior cycles. 

An alternative design \citep{res2010} is to update the elite set dynamically.  
The solutions obtained with the application of PR are considered candidates
to update the elite set. PR is applied as long as there exists a pair of elite solutions that have not yet been relinked, as shown
in step \ref{alg:alg-4:while} of Algorithm~\ref{alg:alg-4}. 

\begin{algorithm}[htbp]
	\caption{\texttt{GRASP with Evolutionary PR}}
	\label{alg:alg-4}
	\begin{algorithmic}[1]
		\State{Apply dynamic GRASP with PR} 
		\State{Let $\mathit{ES}$ be the resulting elite set}
		\While {$\exists \; S,T \in \mathit{ES}$ not yet relinked} \label{alg:alg-4:while}
		\State{$S \gets \texttt{PR}(S, T)$} 
		 \State{$ \mathit{ES} \gets \texttt{Update\_Elite}(ES,S) $}
		\If{$f(S) > f_*$} \label{alg:alg-4:ifPR}
		\State{$ S_* \gets S $} \label{alg:alg-4:updateSs2}
		\State{$ f_* \gets f(S_*) $} \label{alg:alg-4:updatefs2}
		\EndIf
		\EndWhile 
		\State \Return $S_*$ \label{alg:alg-4:return}
	\end{algorithmic}
\end{algorithm}

\section{Historical evolution and critical review}
\label{sec:history}

As aforementioned, GRASP with Path Relinking was first proposed by \cite{lag1999} to minimize arc crossings in 2-layered networks. In this problem, a solution is represented by two permutations, which determine the ordering of the vertices in the two layers of a hierarchical drawing. The absence of additional constraints makes the application of Path Relinking particularly suitable for this problem. However, it is important to note that when dealing with problems that involve numerous constraints, the application of Path Relinking may become complex. This complexity arises from the need to incorporate additional search mechanisms to ensure feasibility in the intermediate solutions.

The relinking in the original implementation of GRASP with PR involves finding a path between a solution obtained after an improvement phase of GRASP and a randomly selected elite solution. This particular variant is now known as greedy forward PR with a dynamic implementation. For the sake of simplicity, we will refer to the original implementation as OGPR.

Although the authors acknowledge that GRASP solutions are generated independently and initially unrelated (i.e., not linked by a path such as in the case of tabu search), they chose to retain the name of the method, Path Relinking, to give credit to its roots in tabu search. This seminal paper presents two notable contributions. Firstly, the hybridization of two different methodologies: a memory-less approach (GRASP) and a memory-based approach (PR). Secondly, the introduction of an elite set to facilitate communication between them.

At the Metaheuristic International Conference (MIC) held in Porto, Portugal, in 2001, \cite{Binato2001} introduced Greedy Randomized Adaptive Path Relinking (GRAPR), which extends the GRASP principles to PR. Recognizing the existence of an exponentially large number of paths between two solutions, this method samples a set of high-quality paths as discussed in Section~\ref{grpr}.

The annotated bibliography of GRASP literature from 1989 to 2001, documented by \citet{Festa2002}, acknowledged the seminal paper by \cite{lag1999} as a significant enhancement to the basic GRASP methodology. This bibliography, combined with the presentation mentioned at the MIC conference, sparked the interest of numerous researchers in applying PR post-processing to their GRASP methods. 
In this section, we review some of the most important contributions of this new hybrid methodology.
We divide this period into two decades, the creative decade (2000-2010) and the consolidation decade (2010-2020).

\subsection{The creative decade (2000 - 2010)}

During the early years of the decade from 2000 to 2010, many researchers dedicated their efforts to improving the GRASP algorithm by incorporating into it different Path Relinking designs. This integration proved to be effective in addressing a wide range of combinatorial optimization problems as shown next.

Notably, \cite{Ribeiro2002} focused their work on the Steiner problem in graphs and introduced an adaptive Path Relinking approach to enhance their GRASP implementation with the dynamic strategy (refer to Section \ref{dynamic}). The authors considered two distinct types of moves for executing the Path Relinking procedure. Following the application of these moves for several iterations, the adaptive strategy involved selecting the superior move between the two based on its performance within the specific problem instance. This method demonstrated robustness across a diverse range of instances.

\cite{Festa2002b} adapted the original design OGPR to tackle the max-cut problem, a classic NP-hard combinatorial optimization problem. An important contribution of this paper in terms of methodology is that PR is applied as a post-processing procedure to enhance both GRASP and VNS. This shows the effectiveness of coupling different metaheuristics with PR and highlights its potential to enhance them. Although the original hybridization limited itself to couple PR with GRASP, the proposal can be applied to any other metaheuristic, as will be shown in the following years.

It is very interesting to observe that in these early years of GRASP with PR, each new paper includes a methodological contribution, extending the original design to improve the search in the solution space. This is also the case of \cite{aiex2003}, who consider not only quality, but also diversity when updating the set of elite solutions. In particular, if $S_b$ and $S_w$ are the best and the worst solutions respectively in the elite set, $\mathrm{ES}$, obtained with the application of GRASP, a new solution $S$ obtained with PR qualifies to enter in $\mathrm{ES}$, replacing one of its solutions, if one of these two conditions is satisfied:

\begin{enumerate}
  \item $S$ is better than $S_b$.
  \item $S$ is better than $S_w$, and sufficiently different to all the solutions in $\mathrm{ES}$.
\end{enumerate}

If $S$ is accepted to enter the $\mathrm{ES}$, it replaces $S_w$, which will be removed from the set. This implementation requires the definition of a distance between pairs of solutions, and a threshold to determine whether two solutions are similar. \citet{aiex2003} also explored parallel designs on multiprocessors of GRASP with PR, with the classical collaborative and non-collaborative implementation that here translates into information exchange among processors regarding the solutions in the elite set. This is elaborated in \citet{Aiex2005b} and coined as \emph{independent GRASP with PR} and \emph{cooperative GRASP with PR}. 

Despite being, at that time, a recent method, \citet{Resende2003} suggest certain strategies, based on their experience, to efficiently apply GRASP with PR, and test these strategies in the context of virtual circuit routing in private networks. In particular the authors recommend:

\begin{itemize}
    \item Not to apply PR at every GRASP iteration, but periodically;
    \item Given two solutions, explore only one trajectory to connect them;
    \item Do not follow the full trajectory, but instead only examine it partially.
\end{itemize}

The authors were actually setting the foundations of what we call now \emph{forward} and \emph{backward} PR as well as \emph{truncated} PR (see Section \ref{sec:AdvancedPR}). Specifically, they observed that better solutions are found when PR starts from the best of the two solutions to be relinked given that PR more carefully explores the neighborhood of the best solution.

In their study, \cite{aiex2005} applied the original dynamic PR implementation, OGPR, proposed by \cite{lag1999}, to tackle the three-index assignment problem, an extension of the classical two-dimensional assignment problem to a tripartite graph. The authors followed the standard design where each solution generated with GRASP is relinked with an elite solution. Moreover, they introduced two innovative strategies into their PR implementation. The first strategy employed what now is recognized as forward and backward PR. Specifically, for each pair of solutions, it considered two paths, one in each direction.
The second strategy is an application of static PR,
involved a post-processing step that effectively relinked all pairs of solutions within the elite set. Notably, what stands out is that rather than selecting between a dynamic or static implementation, this approach interleaves both strategies.

The variant called Evolutionary Path Relinking (EvoPR) was proposed by \cite{res2004}. As described in their paper, the motivation to create this variant was to augment PR with the concept of multiple generations taken from Genetic Algorithms (GAs). Although the method was not coined in that paper, the algorithm to iterate over the set of elite solutions, was proposed there. Interestingly, the evolutionary design proposed does not share with GAs the typical sampling mechanism based on randomization; instead of that, the authors proposed an exhaustive design in which all the pairs of solutions in the elite set are subject to PR in exactly the same way that combinations are performed in Scatter Search (SS) \citep{lag2003}. As a matter of fact, EvoPR and SS share the same design in terms of pool management of the elite set. From a taxonomy of heuristics perspective, we may consider EvoPR as a particular case of SS, since in the latter we can apply PR or any other combination method to generate new solutions. It is worth mentioning however, that PR is usually applied in SS to solve combinatorial optimization problems, while other combination methods may be used in SS when the problem on hand is a nonlinear (global) optimization problem.

In 2005, we may find the first paper summarizing the alternatives and strategies proposed so far to hybridize GRASP with PR. \cite{ResendeRibeiro2005} outline their tutorial talk on GRASP with PR, given at the Metaheuristics International Conference (MIC 2003) in Kyoto. In this paper they first summarize the different ways in which memory-based mechanisms have been included in GRASP to enhance its performance with a long-term memory structure: 
\begin{itemize}
    \item 
    Filtering low-quality constructions to skip the application of the local search;
    \item 
    Checking for duplicated constructions with a hash table;
    \item 
    Identifying the strongly-determined variable values in high-quality (elite) solutions to construct new good ones;
    \item 
    Dynamically tuning the $\alpha$ parameter in GRASP according to its performance in past iterations;
    \item 
    Apply PR to elite solutions.
\end{itemize}      
Then, the authors chronicle seven implementations of PR: periodical, forward, backward, back and forward, mixed, randomized, and truncated. 
Mixed path relinking was initially proposed by \citet{Glover1996b}. Its first application and testing were in the setting of the 2-path network design problem by \citet{Rosseti2003}, followed by other studies \citep{ResendeRibeiro2005,RibRos2009a}. These studies demonstrated that mixed path relinking generally surpasses other methods such as forward, backward, and back-and-forward path relinking in terms of the trade-off
between runtime and solution quality.
We refer the reader to Section \ref{sec:AdvancedPR} for a detailed description of these strategies. Finally, Resende and Ribeiro describe the two ways in which GRASP can be hybridized with PR through an elite set, as an intensification mechanism or as a post-processing after the termination of GRASP:

\begin{itemize}
    \item Intensification PR: It is applied to link a generated (and improved) GRASP solution with an elite solution, as in the original design. PR is applied in each GRASP iteration step; 
    \item Post-optimization PR: It is applied to all pairs of solutions collected in an elite set obtained with the application of GRASP. Here PR is run once GRASP has been applied (and finished).
\end{itemize}

\cite{ResendeRibeiro2005} perform some experiments and conclude that the best results are obtained when we combine both strategies, intensification and post-optimization PR, which actually constitutes the design of Evolutionary PR. It should be noted that these good results are obtained at the expenses of a computational cost and long running times. Their findings are in line with other PR implementations. For example, \cite{Pinana2004} applied a dynamic PR post-processing to GRASP for the matrix bandwidth minimization with excellent results. Their method is actually Evolutionary PR although at that point in time the term had not been popularized yet.

The good results obtained with Path Relinking when applied to elite GRASP solutions, encouraged many researchers to try it with their previous GRASP algorithms, improving their best-known heuristic results. This is the case for example of \cite{Festa2007} for the weighted Max-SAT problem, that was previously solved by GRASP. In these studies, a key question was to prove if it is worthwhile to invest some extra time applying PR after GRASP or it is better to simply apply GRASP for the total amount of time. This question was answered affirmatively in many papers like this one, concluding that at some point GRASP stagnates and it is worth to apply then PR to reach higher quality solutions.

\subsection{The consolidation decade (2010 - 2020)}

During the first decade in our revision, from 2000 to 2010, many presentations in heuristic conferences were increasingly devoted to path relinking hybridizations, not only with GRASP but also with other metaheuristics.  It is worth mentioning that the relevance of GRASP helped Path Relinking to gain prominence and become a usual post-processing for many other metaheuristics since then. For example, the well-known Genetic Algorithms also benefit from the addition of Path Relinking, sometimes as a combination method, as in \cite{VALLADA2010} for the flowshop problem.

At the end of the first decade in our study (2000 - 2010), most of the designs described in the previous section were clearly established and researchers in the second decade, from 2010 to 2020, mainly applied these designs to different optimization problems. We may say that the paper by \cite{res2010} on the max-min dispersion problem (MMDP), is indeed a review paper in which some of these previous designs are labeled and tested, and can be considered a milestone in the historical evolution of the methodology. In particular, the authors tested four
different variants of PR: greedy PR, greedy randomized PR,
truncated PR, and evolutionary PR, as well as two search strategies: static and dynamic. Their experiments showed that the dynamic variants of GRASP with greedy PR and
GRASP with evolutionary PR are the best methods for the
MMDP instances tested. Moreover, the
results indicate that the proposed hybrid heuristics compare
favorably to previous metaheuristic algorithms for the MMDP, such as tabu search and simulated annealing. 

\cite{res2010} pointed out that the results obtained with their implementation were not all due to the strategies tested, but were also dependent (and 
enhanced) by the context-specific methods developed for the MMDP.  Measuring the contribution of a strategy for a specific implementation for a given problem is a cornerstone in heuristic optimization.  

Most of the GRASP with PR publications in the second decade limit themselves to a parti\-cular design already proposed, with a few exceptions that we describe now. Among the numerous papers we highlight the following ones that propose if not new, somehow generic PR strategies or implementations.

\begin{figure}[ht]
\centering\includegraphics[angle=0,scale=0.75]{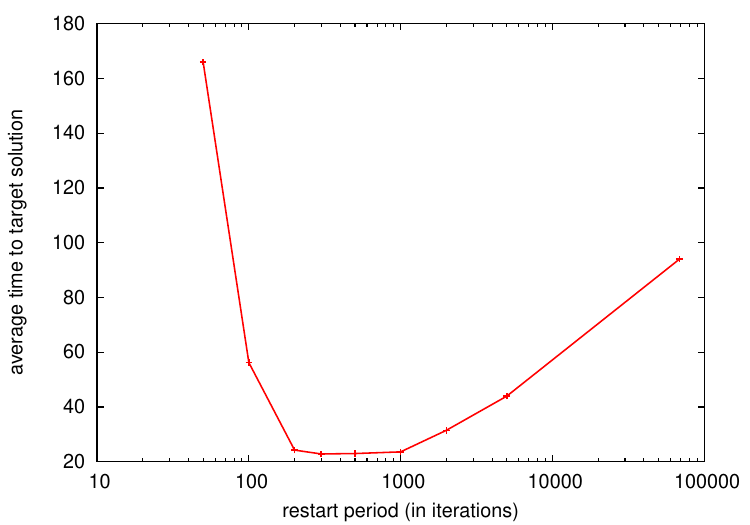}
\caption{
Average time to target solution using different
values of $\kappa$. }
\label{f_g12-avg-iter-ttt}
\end{figure}

The re-start strategies proposed by \cite{Resende2011} deserve special mention. These strategies are used to re-initiate the search from scratch. In short, the strategy can be described as: Keep track of the last iteration when the incumbent solution was improved and re-start the GRASP with
path-relinking heuristic if $\kappa$ iterations have gone by without improvement. Re-starting GRASP with path-relinking requires emptying out the elite set, and starting a new iteration with a new random number generator seed. While this strategy also requires us to input a value for parameter $\kappa$, a limited number of values for $\kappa$ almost always achieves the desired result, i.e. to reduce the average iteration count as well as its standard deviation.
Figure~\ref{f_g12-avg-iter-ttt} from \cite{Resende2011} shows
average time to target solution for maximum
cut instance $\mathit{G12}$ using different
values of $\kappa$ in GRASP with PR with restart. All runs of all strategies have found a
solution at least as good as the target value of 554.

An important contribution in this decade is the analysis of the different adaptations of GRASP with PR to multiobjective optimization \citep{mar2015}. The authors first reviewed previous applications of both metaheuristics to then establish a classification of the different ways in which they can be applied. The first strategy, \textit{pure}, is one where each objective is optimized in isolation. In a second strategy, \textit{sequential}, each objective alternates to guide the search.  The last strategy, \textit{weighted}, is one where all the objectives are combined into a single
master objective. The practical side of the study was carried out  with two hard biobjective combinatorial problems, the path dissimilarity (PDP) and the bi-orienteering problems (BOP), to test the 26 different GRASP and PR variants proposed in the paper. An interesting conclusion of the study is that in each problem the best results are obtained with a different GRASP with PR variant. Specifically, for the PDP, the pure variant achieves the best results, while in the BOP the weighted variant is the winner.

Exterior PR has been largely overlooked in spite of the large number of papers applying PR. Two exceptions are \cite{duarte2015greedy} and \cite{rodriguez2017grasp}, in which it is  applied to achieve additional diversification.  As a matter of fact, the experiments in these two papers indicate that Exterior PR is indeed a powerful strategy to achieve diversification, while at the same time obtaining good solutions (in terms of the objective function). More recently, \cite{loz2023b} propose a static GRASP that incorporates a hybrid approach combining Interior and Exterior PR for a bi-objective problem. In particular, if the initiating and guiding solutions are similar, Interior PR is applied to intensifying the search; otherwise, Exterior PR is applied to enhance search diversification. 

Stochastic optimization has also been approached with GRASP with Path Relinking. \cite{albareda2017heuristic} adapt it to solve the faciliy location problem with stochastic demands. GRASP is applied two create two pool of solutions. The elite pool with the best solutions found with GRASP, and the diverse pool, which consists of solutions that were not good enough to be part of the elite pool but that are the most diverse ones found. The path relinking is applied from an elite solution to a solution randomly selected from the diverse pool. In our view, PR can deal very well with the stochastic nature of uncertain data, and we will probably witness more applications like this one in the following years.

\begin{figure}[ht]
\centering
\includegraphics[width=0.6\textwidth]{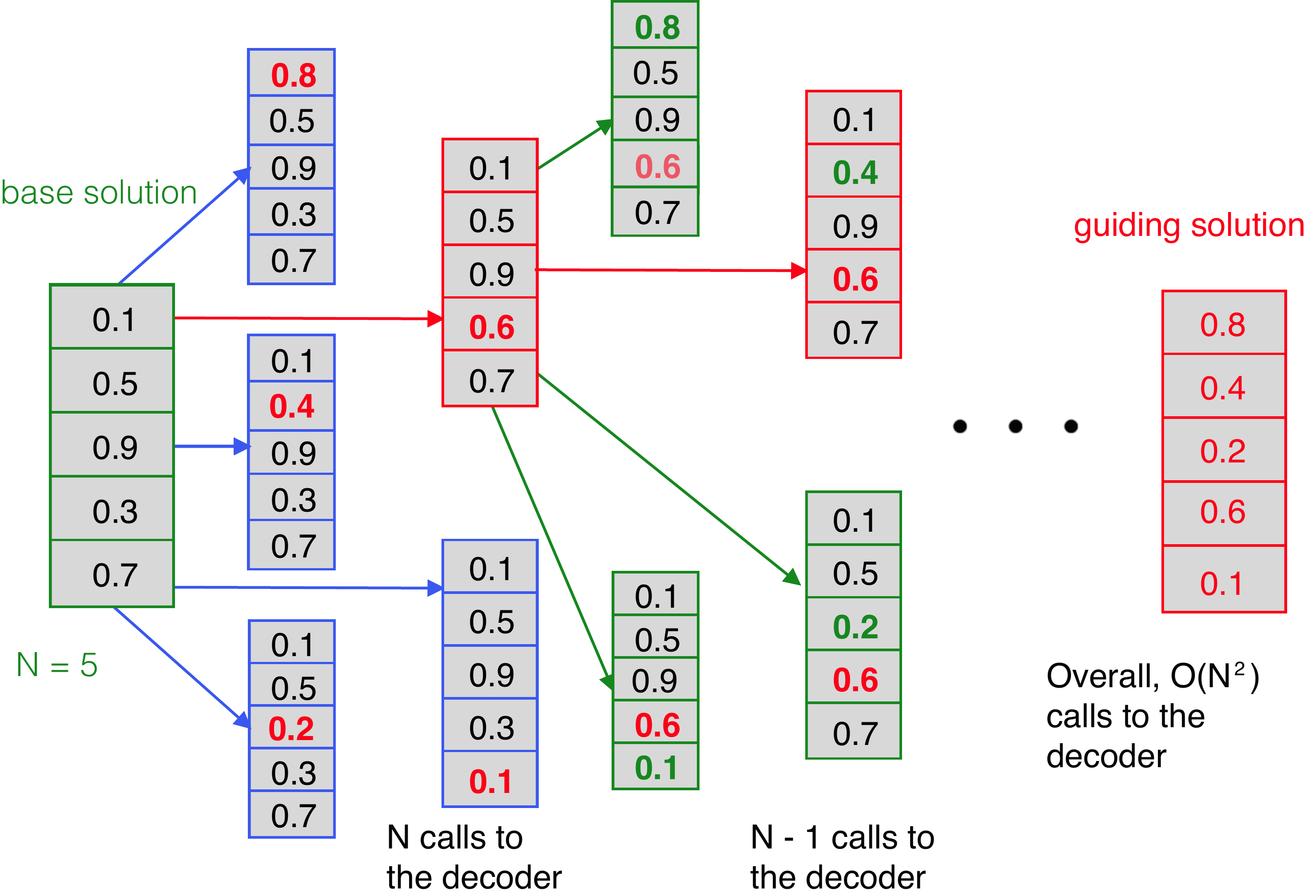}
\caption{
Implicit path relinking. }
\label{f_ipr}
\end{figure}

Path relinking has found its way into Biased Random-Key Genetic Algorithms (BRKGA) for combinatorial optimization in \citet{AndTosGonRes21a}.  A BRKGA works in the half open unit hypercube where each point in the
hypercube is an encoded solution which is decoded into the solution space of the combinatorial optimization problem by applying a decoder \citep{GonRes11a}. \textit{Implicit Path Relinking}, or IPR, applies the concept of Path Relinking to link two solutions encoded by random-key vectors in the half open unit hypercube in $\mathbb{R}^n$.
In its simplest form, IPR has an initiating solution $S$ and a guiding solution $T$, each a vector of random keys. These could be, for example, two elite solutions in distinct islands in an island model BRKGA. To begin, the base solution $S'$ is set to the initiating solution, i.e. $S' \gets S$.
At each step of IPR, a random key of the guiding solution replaces the corresponding key in the base solution. To evaluate all key insertions requires $n$ calls to the decode in the first iteration, $n-1$ calls in the second iteration, and so on. The best-evaluated solution is accepted as the next base solution and the procedure continues with the next iteration. Not only is this approach slow, requiring $O(n^2)$ calls to the decoder, it focuses excessively on intensification, not having enough diversification. A hidden structure should be expected to appear in the vectors of random
keys throughout the evolutionary process. 
Figure~\ref{f_ipr} illustrates two iterations of IPR. In the first iteration each weight
of the guiding solution replaces the corresponding weight in the base solution.  The best-evaluated solution, $S'= (0.1, 0.5, 0.9, 0.6, 0.7)$, is indicated by the incoming red arrow. For the second iteration, $S'$ is made the base solution where the fourth weight is now fixed at value 0.6 and weights of the guiding solution replace unfixed weights (shown in green). The best-evaluated solution, $S''= (0.1, 0.4, 0.9, 0.6, 0.7)$, is indicated by the incoming red arrow. This process continues until the guiding solution is reached.

Certain portions of the random-key vectors can become essential building
blocks for high-quality solutions. Tearing these blocks apart when combining with other solutions, can destroy
the high-quality solutions. \citet{AndTosGonRes21a} propose an approach that exchanges a block of keys (rather than a single
key).
The number of blocks $\mathit{NB}$ depends on the the block size and the length of the
path. In this approach, the first iteration requires $\mathit{NB}$ calls to the decoder, the second iteration requires $\mathit{NB}-1$ calls, and so on. The total number of calls is $O((\mathit{NB})^2)$.
IPR has be applied to other random-key optimizers, such as the ones described in \citet{rko-robot-2022}.

\cite{Sanc2021} propose a GRASP with PR for multi-objective problems that incorporates strategic oscillation (SO) to enhance the exploration of different regions in the solution space. Specifically, the algorithm solves a multiobjective facility location problem where  $k$ facilities have to be located among $m$ potential ones. Initially, the algorithm generates the set of efficient solutions using GRASP, alternating the objective function in each construction. Subsequently, SO is applied to each efficient solution previously obtained. For each solution and each objective, the algorithm obtains two infeasible solutions by greedily adding and dropping some elements. Then, PR is applied between each pair of solutions as a repairing mechanism to recover feasibility. This implementation combines PR and SO in an efficient way that crosses back and forth the feasibility boundary, and may be generalize to other constrained problems.

An open question after all the research described above is if we can establish which strategies and variants consistently produce the best results or they are problem dependent and cannot be set beforehand. Our computational experience in Section \ref{sec:comput} helps to clarify this question.

\section{Connections with other methodologies}
\label{sec:other_methods}

A general classification of metaheuristics consists of positioning the methods into trajectory-based or population-based. Trajectory-based methods operate on a single solution and employ neighborhood searches to identify a transition (or move) from a current solution to the next. Population-based procedures, on the other hand, operate on a set of solutions throughout the search. GRASP and basic tabu search are two examples of trajectory-based approaches. However, when path relinking is added to these methods, the classification is not as obvious, given that PR requires a set of elite solution that need to be maintained during the search.

GRASP with PR has strong connections with the population-based method known as scatter search (SS) that uses path relinking to combine solutions. Here are the main connections:

\begin{itemize}
   \item SS starts with a population of solutions that is constructed to balance quality and diversity. Many SS implementations use GRASP constructions for this purpose. This initial population is equivalent to the solutions found in the construction phase of a GRASP with PR.  
   \item SS maintains a reference set of solutions that may be viewed as being equivalent to the elite set used in GRASP with PR.
   \item The combination method in SS is a path relinking of reference solutions. Since only new reference solutions are combined with others, this is equivalent to combining a new local optimum found after a GRASP iteration with an elite solution.
\end{itemize}

Despite these close connections, a SS that uses PR as its combination method and a GRASP with PR are not expected to perform identically. One major difference comes from the role of GRASP constructions in each method. Suppose that GRASP constructions are used within scatter search to create an initial population of 100 solutions. By SS rules, only a small fraction (e.g., 10\%) of these solutions are chosen to build the initial reference set. The remaining solutions in the population are engaged only if the SS includes a re-starting mechanism that is called when the reference set is rebuilt after it converges. In contrast, a GRASP with PR that performs 100 iterations employs all the constructions (one by one) to find local optima that are turned into initiating solutions in the PR process. Additional differences are associated with the rules for updating the reference set in SS versus the updating of the elite set in GRASP with PR.

Connections between GRASP with PR and tabu search with PR can also be established. For instance, some tabu searches are designed to start from a GRASP construction. This is followed by a sequence of neighborhood searches that ultimately reach a local optimum. This first local optimal point is added to the elite set that is being built and maintained for path relinking. Short-term memory strategies are then used to escape local optimality and move the search to new basins of attraction. As this process continues, the set of elite solutions is populated and updated. Path relinking in the context of tabu search has been implemented in various ways:

\begin{itemize}
   \item A PR phase is triggered every time the search reaches a local optimum, where the local optimal solution becomes the initiating solution and the guiding solution is chosen from the elite set (see e.g., \cite{Arment11})
   \item A PR phase is called within the tabu search every fixed number of iterations (see, e.g., \cite{dib2018}
   \item PR is a post-processing phase that operates on a population of local optimal solutions found with tabu search (see e.g., \cite{peng2015}) 
\end{itemize}

The first design from the list above is the one with the strongest connection to the original GRASP with PR proposal by \cite{lag1999}.  
Since PR can in theory be applied to any metaheuristic algorithm that generates a sequence of locally optimum solutions, its hybridization with such algorithms should be further explored since we have observed little development in this direction, as compared with GRASP.

\section{Numerical experiments}
\label{sec:comput}

In this section we present experimental results with several strategies of GRASP and GRASP with
Path-Relinking from Sections ~\ref{sec:origins}, \ref{sec:AdvancedPR} and~\ref{sec:GRASP_PR}:
\begin{itemize}
\item Multi-Start Semi-Greedy algorithm; 
\item GRASP with first and best improvement local search;
\item Forward, Backward, and Mixed PR;
\item Static, Dynamic, and Evolutionary PR;
\end{itemize}
Additionally we tested the application of 
Local Search to all path solutions, only to the best path solution, and to periodically determined solutions, e.g. every $q$ solutions in the path.

Each algorithm is run on each instance using two time horizons, 60 and 600 seconds on an AMD Epyc 7282 2.8 GHz processor with 32 GB of RAM. All algorithms were coded in Java 17.  We perform seven experiments and summarize the results for each experiment in the tables below. In each table, we report average values for the following statistics: 
\begin{itemize}
     \item \textbf{$\#Best$}: Number of instances in which the method obtains the best solutions in each individual experiment;

     \item \textbf{$\% Dev$}: Average percentage relative deviation between the  value obtained with the method on each instance and the best value obtained for this instance in this experiment;
  
    \item \textbf{$\#Best_k$}: Number of instances in which this method obtains the best known solution across all experiments;

    \item \textbf{$\% Dev_k$}: Average percentage relative deviation between the  value obtained with this method on each instance and the best known value for this instance across all experiments.

\end{itemize}

\subsection{Problems and instances}

We consider two well-known NP-hard problems in combinatorial optimization to perform our experiments with the proposed GRASP with Path Relinking variants: The linear ordering problem (LOP) and the max-cut problem (MAX-CUT). 

In its graph version, the LOP is defined as follows.  Let $G = (V,A)$ denote a complete digraph with node
set $V=\{1,2,\ldots,n\}$ and arc set $A$ containing for every pair of nodes $i$ and $j$ arcs $(i,j)$ and $(j,i)$. Let $c_{ij}$ be the cost of arc $(i,j)$. Note that $c_{ij}$ is usually not equal to $c_{ji}$. 
A tournament~$\tau$ in $A$ consists of a subset of arcs containing for every pair of nodes $i$ and $j$ either arc
$(i,j)$ or arc $(j,i)$, but not both. A \emph{spanning acyclic tournament} is a tournament without directed cycles, i.e., not containing an arc set of the form $\{(v_1,v_2),(v_2,v_3),\ldots,(v_k,v_1)\}$ for some $k>1$
and distinct nodes $v_1,v_2,\ldots,v_k$.
The LOP can be stated in terms of graphs as computing a spanning acyclic tournament $\tau$ in $A$ such that $\sum_{(i,j)\in \tau}c_{ij}$ is as large as possible.

Alternatively, the LOP can also be defined on a matrix,
the so-called \textit{triangulation problem}. Given an $n \times n$ matrix $C=(c_{ij})$ determine a
simultaneous permutation of the rows and columns
of~$C$ such that the sum of the upper triangular entries, $\sum_{i<j} c_{ij}$, becomes as
large as possible.

A linear ordering, or permutation, of the nodes $\{1,2,\ldots,n\}$ is a feasible solution of the LOP. 
In \cite{marti2011}, it is shown that by setting arc weights for the complete
digraph equal to the matrix entries, the triangulation problem and the spanning acyclic tournament problem are equivalent. The LOP has applications in aggregation of individual preferences, ranking in sports tournaments, and scheduling with precedence among others.

\cite{marti12} compiled a benchmark library for the LOP, the so-called LOLIB with 484 instances. We consider in our study 80 instances from two widely used classes of instances in the LOLIB, the 50 input-output I/O tables (matrices) generated from economic sectors in different regions, and the 30 MB instances \citep{mit2000}. 
The MB instances are random matrices where the lower-diagonal entries are uniformly distributed in [0, 99] and the upper-diagonal entries are drawn uniformly from [0, 39].

The MAX-CUT problem can be described on an undirected graph $G=(V,E)$, where $V=\{1,\dots,n\}$ is the set
of vertices and $E$ is the set of edges,
and weights $w_{ij}$ 
associated with the edges $(i,j)\in E$. The MAX-CUT problem
consists of finding a 
subset of vertices $S$ such that the weight of the cut $(S,\bar S)$ given
by $$w(S,\bar S)=\displaystyle{\sum_{i\in S,j\in \bar{S}} w_{ij}}$$
is maximized, where $\bar S=V \setminus S$.

Applications are found in VLSI design and statistical physics, see e.g. 
\citet{BaGrJuRe88,ChaDu87,CheKajCha83,Pint84}, among others.
The reader is referred to \citet{PolTuz95} for an introductory 
survey. 

For the experiments with MAX-CUT, we used three sets of instances used in \citet{Festa2002b}. The first set of 54 MAX-CUT instances were created
by \citet{HelRen97}.
The second set of instances is based on instance pm3-8-50, from
the benchmark problem set of the 7th DIMACS Implementation Challenge\footnote{\url{http://archive.dimacs.rutgers.edu/Challenges/Seventh/}}. These four
instances were generated by M. J\"{u}nger
and F. Liers using the Ising model of spin glasses. 
Finally, the last set of 30 instances
arise in physics, and were
proposed by \citet{BurMonZha01}.

Note that a solution to the LOP is a permutation whereas a solution to the MAX-CUT is a partition represented as a binary vector. Our experiments therefore cover these two types of combinatorial optimization types.

\subsection{GRASP experiments} \label{sec:GraspExp}

This experiment focus on the traditional GRASP implementation. Its purpose is to empirically identify the best components for both the construction and improvement phases.

\vspace{0.2cm} 
\underline{Experiment 1}. At each step of the construction phase of GRASP, an element is randomly selected from as so-called restricted candidate list (RCL) to be added to a partial solution. Constructions are greedy or semi-greedy because membership to RCL is determined by a greedy function. Considering that both of test problems can be represented in a graph $G$, let $g(v)$ represent the ``attractiveness" of adding vertex $v$ to a partial solution being constructed. In this experiment, we test two procedures to build the RCL. Both procedures are based on first creating a candidate list (CL) of elements (in our case vertices). The candidate list starts with all vertices in the graph, that is $CL=V$. The solution $S$ under construction starts empty. For the LOP, $S$ is an ordered list of vertices and a feasible solution is reached when $|S|=n$, that is, when all vertices have been added. For the MAX-CUT problem, $S$ is a subset of vertices, and any subset represents a feasible solution. Since the LOP and the MAX-CUT are maximization problems, we order the candidate list by decreasing $g(v)$ value, where $g(v)$ is the increase in the objective function value $f(S)$.

The first is a threshold-based procedure (SemiGreedy1) that uses a parameter $\alpha$ to control the size of the RCL. Let $g_{max}=\max\limits_{v \in CL} (g(v))$ be the maximum attractiveness value, that is, the one corresponding to the first vertex in CL. Then, the RCL consists of all vertices in CL whose attractiveness is $\alpha\%$ away from the maximum attractiveness:

$$RCL=\{v \in CL:g(v) \ge (1 - \alpha) g_{max}\}$$

When $\alpha = 0$ in this RCL, the solution construction is totally greedy and when $\alpha = 1$ the construction is totally random. Instead of a fixed $\alpha$ value, our procedure, at each step of the construction process, randomly selects one value in the range $(0,0.3]$. That is, we do not consider selecting vertices whose attractiveness is more than 30\% away from the maximum available attractiveness. After the vertex is selected, it is added to $S$ and deleted from CL. Note that for the LOP, the chosen vertex is added at the end of the ordered list.

The second is a cardinality-based procedure (SemiGreedy2) for which $\alpha$ represents the fraction of the most attractive vertices in CL that are included in RCL. Let $p(v)$ be the position of vertex $v$ in CL, where $p(v)=1$ means that $g_{max}=g(v)$. Then, we use the $\alpha$ value to calculate the number of vertices $p_{max}=1+\lfloor \alpha (|CL|-1) \rfloor$ to be included in RCL:

$$RCL=\{v \in CL:p(v) \le p_{max}\}$$

Here again, when $\alpha = 0$ in this RCL, the solution construction is totally greedy and when $\alpha = 1$ the construction is totally random. And, as above, we do not fix the $\alpha$ value and instead we randmly select one value in the range of $(0,0.3]$ at each step of the construction process.

Table \ref{tab:constructMethods} shows the results associated with running these construction procedures for 60 and 600 seconds. The results are in line with the sampling nature of GRASP constructions. That is, as the number of samples (solution constructions) increases, the quality of the best construction increases. The quality metrics indicate that there is no significant difference between SemiGreedy1 and SemiGreedy2 for the LOP. For the MAX-CUT, SemiGreedy1 is the preferred procedure when focusing on solution quality. These two conclusions have been verified by the results of Wilcoxon's paired difference test. The metrics indicate that the construction process is more effective in producing solutions of a higher quality for the LOP than for the MAX-CUT. LOP solutions are represented as permutations while MAX-CUT solutions are represented by binary variables. This may help explain the difference in the performance of the construction procedures. However, the most influential factor can be attributed to the differences in the greedy functions used for the LOP and the MAX-CUT. In both cases, however, SemiGreedy1 performs slightly better than SemiGreedy2 and therefore we choose SemiGreedy1 as our construction procedure.


\begin{table}[htb]
    \centering
    \caption{GRASP constructions.}
\label{tab:constructMethods}
\begin{tabular}{l c c c c c c c c c}
    \toprule
    & \multicolumn{4}{c}{$60$ s. } & &  \multicolumn{4}{c}{$600$ s. } \\
      \cline{2-5}  \cline{7-10} 
    \bf Procedure & \multicolumn{1}{c}{$\# Best$} & \multicolumn{1}{c}{$\% Dev$} & \multicolumn{1}{c}{$\# Best_k$} & \multicolumn{1}{c}{$\% Dev_k$}  & & \multicolumn{1}{c}{$\# Best$} & \multicolumn{1}{c}{$\% Dev$} & \multicolumn{1}{c}{$\# Best_k$} & \multicolumn{1}{c}{$\% Dev_k$} \\
   \midrule
   &  \multicolumn{9}{c}{LOP}\\
   \cline{2-10} 
    {\tt SemiGreedy1} &  49& 0.014&  16& 0.160& & 54& 0.012& 26& 0.120\\
    {\tt SemiGreedy2}  & 45& 0.016& 14& 0.162& & 49& 0.011& 25& 0.118\\
    \midrule
    
    & \multicolumn{9}{c}{MAX-CUT}\\
    \cline{2-10} 
    {\tt SemiGreedy1} &  \multicolumn{1}{c}{88}&  0.000 &  \;3& 7.782& & 88& 0.000&  \;7& 7.226\\
    {\tt SemiGreedy2}  &  \multicolumn{1}{c}{\;1}&  7.433&  \;0 & 14.428& & \;0& 10.168& \;0& 16.414\\

    \bottomrule
\end{tabular}

\end{table}





\vspace{0.2cm} 
\underline{Experiment 2}. The second phase of GRASP focuses on improving the solution that has been built in the first phase. There are two main design decisions to configure a neighborhood search for this phase. First, we need to select a move to transition from one solution to the next. We considered two moves, swap and insert. A swap in the LOP exchanges the position of two vertices in the ordered list $S$. A swap in the MAX-CUT exchanges a vertex in $S$ with a vertex in $\bar S$. In preliminary testing, we determined that this move consistently produced inferior results than a move based on insertions. An insert in the LOP is a move that transfers a vertex from its current position to another position in the ordered list $S$. An insert in the MAX-CUT transfers a vertex in $S$ to $\bar S$ or makes a transfer in the opposite direction.

Establishing the depth of the search is the second decision associated with configuring a neighborhood search. The deepest search, so-called best improving, is to evaluate all moves and selecting the one with the maximum evaluation. This guarantees the steepest accent in a maximization problem. Alternatively, the first-improving strategy stops the exploration as soon as a move that leads to a better solution than the current solution is identified. For the LOP, we stop as soon as an insertion is evaluated to improve the current solution. For the MAX-CUT, we implemented the first-improving strategy suggested by \cite{mardualag09} and that they refer to as $LS_{2}$.

Table \ref{tab:grasp} shows the results of applying the full GRASP to the LOP and MAX-CUT instances. As the table shows, there is no detectable differences in performance between the best-improving and the first-improving strategies when tackling LOP instances. For the MAX-CUT, there is a significant difference between  best-improving and the first-improving when considering relative deviation of the solutions against the best known values. We tested for statistical significance in both cases using Wilcoxon's paired difference test. In light of these results, we configured our GRASP with SemiGreedy1 constructions, insert moves, and best-improving exploration depth. 

We also performed a one-tailed test to detect significant differences between the results obtained in Table \ref{tab:constructMethods} and the results reported in Table \ref{tab:grasp}. The purpose of this test was to show that the local search phase is capable of significantly improve upon the solutions created by the construction process. The test rejected the null hypothesis that the mean objective function value of the solutions in Table \ref{tab:grasp} are less than or equal to the mean objective function value of the solutions in  Table \ref{tab:constructMethods}.

\begin{table}[htb]
    \centering
    \caption{GRASP constructions with local search.}
\label{tab:grasp}
\begin{tabular}{l c r r r r r r r r}
    \toprule
    & \multicolumn{4}{c}{$60$ s. } & &  \multicolumn{4}{c}{$600$ s. } \\
      \cline{2-5}  \cline{7-10} 
    \bf Procedure & \multicolumn{1}{c}{$\# Best$} & \multicolumn{1}{c}{$\% Dev$} & \multicolumn{1}{c}{$\# Best_k$} & \multicolumn{1}{c}{$\% Dev_k$}  & & \multicolumn{1}{c}{$\# Best$} & \multicolumn{1}{c}{$\% Dev$} & \multicolumn{1}{c}{$\# Best_k$} & \multicolumn{1}{c}{$\% Dev_k$} \\
   \midrule
   &  \multicolumn{9}{c}{LOP}\\
   \cline{2-10} 
    {\tt GRASP (First)} &  55& 0.004&  43& 0.018& & 58& 0.002& 47& 0.011\\
    {\tt GRASP (Best)}  & 69& 0.001& 46& 0.015& & 70& 0.001& 49& 0.009\\
    \midrule
    
    & \multicolumn{9}{c}{MAX-CUT}\\
    \cline{2-10} 
    {\tt GRASP (First)} &  \multicolumn{1}{c}{88}&  0.000 &  14& 0.905& & 88& 0.000&  15& 0.628\\
    {\tt GRASP (Best)}  &  \multicolumn{1}{c}{\,8}&  4.986&  8& 5.826& & 11& 4.649& 11& 5.234\\

    \bottomrule
\end{tabular}

\end{table}


\subsection{PR experiments} \label{sec:PRexp}

The experiments in this section have the goal of detecting performance differences among the various path relinking configurations described above. For this experiment, we use GRASP as a means for generating one hundred solutions. We then use a static path relinking ({\tt StPR}) design in which the top ten solutions (according to quality) are selected from the set of one hundred to form the elite set. We then apply PR to all pairs of solution in the elite set and return the best solution found (see Section \ref{sec:GRASP_PR}). 


\vspace{0.2cm} 
\underline{Experiment 3}. This experiment assesses the performance of the three mechanisms to transition from the current solution in the direction of the guiding solution during the path relinking process. These mechanisms, labeled  {\tt Forward}, {\tt Backward}, and {\tt Mixed} are described in Section \ref{sec:AdvancedPR}. Table \ref{tab:StaticPathRelinking} summarizes the results of this experiment. The statistics shown in the table show small differences in performance among the three designs for the LOP. There is a slight advantage of a {\tt Forward} path relinking in the context of the MAX-CUT problem. However, this advantage is relatively small. Ideally, we would like to provide some insights about this result (i.e., that the three PR mechanisms seem to perform at the same level), but our analysis of the behavior of the algorithm and the paths that it generated did not result in a solid conjecture.




\begin{table}[htb]
    \centering
    \caption{Static path relinking with three search directions.}
\label{tab:StaticPathRelinking}
\begin{tabular}{l c c c c c c c c c}
    \toprule
    & \multicolumn{4}{c}{$60$ s. } & &  \multicolumn{4}{c}{$600$ s. } \\
      \cline{2-5}  \cline{7-10} 
    \bf Procedure & \multicolumn{1}{c}{$\# Best$} & \multicolumn{1}{c}{$\% Dev$} & \multicolumn{1}{c}{$\# Best_k$} & \multicolumn{1}{c}{$\% Dev_k$}  & & \multicolumn{1}{c}{$\# Best$} & \multicolumn{1}{c}{$\% Dev$} & \multicolumn{1}{c}{$\# Best_k$} & \multicolumn{1}{c}{$\% Dev_k$} \\
   \midrule
   &  \multicolumn{9}{c}{LOP}\\
   \cline{2-10} 
    {\tt StPR Forward} &  51& 0.003&  50& 0.007& & 50& 0.002& 50& 0.004\\
    {\tt StPR Backward}  & 67& 0.001& 51& 0.005& & 65& 0.000& 55& 0.002\\
    {\tt StPR Mixed}  & 70& 0.000& 52& 0.004& & 75& 0.000& 55& 0.002\\
    \midrule
    
    & \multicolumn{9}{c}{MAX-CUT}\\
    \cline{2-10} 
    {\tt StPR Forward} &  \multicolumn{1}{c}{68}&  0.040&  17& 0.611& & 62& 0.048&  26& 0.388\\
    {\tt StPR Backward}  &  \multicolumn{1}{c}{43}&  0.122&  15& 0.692& & 46& 0.075& 26& 0.414\\
    {\tt StPR Mixed}  &  \multicolumn{1}{c}{37}&  0.145&  16& 0.715& & 50& 0.136& 29& 0.475\\
    \bottomrule
\end{tabular}

\end{table}

We performed the same experiment using the Dynamic ({\tt DyPR}) and the Evolutionary ({\tt EvPR}) designs. As shown in \ref{append}, our experiments were inconclusive in terms of identifying the best strategy to create a relinking path across both problem types. Using $\#Best$ as the criterion for selecting the best mechanism, the results with {\tt StPR} reported in Table \ref{tab:StaticPathRelinking} would lead us to choose {\tt Mixed} for the LOP and {\tt Forward} for the MAX-CUT. However, these choices change when examining the $\#Best$ values in \ref{append}. What this shows is that there is a significant interaction between the strategies to choose the elite set and the strategies to build the relinking path.



\vspace{0.2cm} 
\underline{Experiment 4}. This experiment tests the use of local search during the path relinking. Recall that an option in PR is to add intensification in the form of applying a local search procedure to intermediate solutions. Applying to all solutions could be computationally wasteful because the difference between two adjacent solutions in the relinking path is a single move. Nonetheless, we tried this alternative (All) along with two other options, one in which the local search is applied to every 5th solution (Every 5), and another where the local search is only applied to the best solution found (Best only). We used {\tt StPR Mixed} for the LOP and {\tt StPR Forward} for the MAX-CUT to produce the results in Table \ref{tab:StaticPathRelinkingLS}.




\begin{table}[htb]
    \centering
    \caption{Static path relinking strategies with local search.}
\label{tab:StaticPathRelinkingLS}
\begin{tabular}{l c c c c c c c c c}
    \toprule
    & \multicolumn{4}{c}{$60$ s. } & &  \multicolumn{4}{c}{$600$ s. } \\
      \cline{2-5}  \cline{7-10} 
    \bf Procedure & \multicolumn{1}{c}{$\# Best$} & \multicolumn{1}{c}{$\% Dev$} & \multicolumn{1}{c}{$\# Best_k$} & \multicolumn{1}{c}{$\% Dev_k$}  & & \multicolumn{1}{c}{$\# Best$} & \multicolumn{1}{c}{$\% Dev$} & \multicolumn{1}{c}{$\# Best_k$} & \multicolumn{1}{c}{$\% Dev_k$} \\
   \midrule
   &  \multicolumn{9}{c}{LOP}\\
   \cline{2-10} 
    {\tt All} &  70& 0.000&  52& 0.004 & & 76& 0.000& 55& 0.002\\
    {\tt Every 5}  & 68& 0.000& 52& 0.004& & 74& 0.000& 54& 0.002\\
    {\tt Best only}  & 70& 0.000& 51& 0.004& & 61& 0.001& 55& 0.003\\
    \midrule
    
    & \multicolumn{9}{c}{MAX-CUT}\\
    \cline{2-10} 
    {\tt All} &  53&  0.085&  17& 0.611& & 48 & 0.093 &  26 & 0.388\\
    {\tt Every 5}  &  80&  0.011&  16& 0.538& & 84 & 0.005 & 31 & 0.300\\
    {\tt Best only}  &  14&  0.377&  14& 0.901& & 16 & 0.376 & 16 & 0.669\\
    \bottomrule
    \multicolumn{3}{l}{0.000 means less than 0.0001}
\end{tabular}

\end{table}

The addition of local search to intermediate solutions does not improve the results for the LOP. It is safe to conclude that, for the LOP, local search should not be applied during the relinking process. It should be applied, if at all, to the best solution only. There is modest improvement in the case of the MAX-CUT when we compare the performance values for {\tt StPR Forward} in Table \ref{tab:StaticPathRelinking} and {\tt Every 5} in Table \ref{tab:StaticPathRelinkingLS} at 600 seconds. In the MAX-CUT case we have verified that there is statistical significance between {\tt Every 5} and {\tt Best only}. Depending on the MAX-CUT context (e.g., when computational time is not an issue), the application of local search to intermediate solutions allowing for some spacing between solution seems to be warranted.


\subsection{GRASP with PR experiments} \label{sec:GRASP-PRExp}

As before, the experiments in this sections are performed on a cumulative design. That is, we use the GRASP configuration with the best performance as we added elements to the design. For the LOP, we use {\tt StPR Mixed Best only} and for the MAX-CUT we use {\tt StPR Forward Every 5}.

\vspace{0.2cm} 
\underline{Experiment 5.} This experiment compares two ways of selecting moves during the relinking process. The selection is made among the changes that make the current intermediate solution move closer to the guiding solution. The first alternative, {\tt GPR}, is to choose, in a deterministic way, the move that results in a solution with the best objective function value ({\tt GPR} stands for Greedy Path Relinking). The second alternative, {\tt GRPR}, is to choose the move randomly among a restricted candidate list, ({\tt GRPR} stands for Greedy Randomized Path Relinking). 




\begin{table}[htb]
    \centering
    \caption{Comparison between deterministic and random move selection.}
\label{tab:MoveSelection}
\begin{tabular}{l c c c c c c c c c}
    \toprule
    & \multicolumn{4}{c}{$60$ s. } & &  \multicolumn{4}{c}{$600$ s. } \\
      \cline{2-5}  \cline{7-10} 
    \bf Procedure & \multicolumn{1}{c}{$\# Best$} & \multicolumn{1}{c}{$\% Dev$} & \multicolumn{1}{c}{$\# Best_k$} & \multicolumn{1}{c}{$\% Dev_k$}  & & \multicolumn{1}{c}{$\# Best$} & \multicolumn{1}{c}{$\% Dev$} & \multicolumn{1}{c}{$\# Best_k$} & \multicolumn{1}{c}{$\% Dev_k$} \\
   \midrule
   &  \multicolumn{9}{c}{LOP}\\
   \cline{2-10} 
    {\tt GPR} &  58& 0.002& 51& 0.004& & 55& 0.002& 55& 0.003\\
    {\tt GRPR}  & 73& 0.000& 60& 0.002& & 80& 0.000& 66& 0.001\\
    \midrule
    & \multicolumn{9}{c}{MAX-CUT}\\
    \cline{2-10} 
    {\tt GPR} &  78& 0.012& 14& 0.905& & 88& 0.000& 31& 0.300\\
    {\tt GRPR}& 88& 0.000& 14& 0.893& & 17& 0.351& 16& 0.649\\
    \bottomrule
\end{tabular}

\end{table}

The results in Table \ref{tab:MoveSelection} once again are inconclusive in terms of an overall strategy for both problem types. Our statistical tests detect performance differences between {\tt GPR} and {\tt GRPR}. However, while {\tt GRPR} is the better strategy for the LOP, {\tt GPR} is the better strategy for the MAX-CUT.

\underline{Experiment 6.} This experiment focuses on GRASP and attempts to detect performance differences among the three strategies for executing PR, namely, {\tt Static}, {\tt Dynamic}, and {\tt Evolutionary}. As previously discussed, the {\tt Static} strategy consists of running GRASP and choosing solutions to create the elite set. Then PR is applied to the elite solutions. In the {\tt Dynamic} strategy, PR is applied after every GRASP iteration, i.e., after the construction and improvement phase result in a local optimal point. The {\tt Evolutionary} variant is also dynamic in nature but the elite set of solutions is updated with solutions that are found during the relinking phase. For the LOP in this experiment we use Greedy Randomized Mixed PR, applying local search to the best solution found in the path. For the MAX-CUT we employ Greedy Forward PR, applying local search to every 5th solution in the relinking path.

Table \ref{tab:StaticDynamicEvol} shows that {\tt Static} and {\tt Evolutionary} have similar performance, while {\tt Dynamic} produces inferior results. This seems to reinforce the merit of working with elite sets that include solutions of the highest quality. This notion is at the core of PR, which as stated by \cite{glolag1997} it is conceived as a ``extreme (highly focused) instance of a strategy that seeks to incorporate attributes of high quality solutions, by creating inducements to favor these attributes in the moves selected."



\begin{table}[htb]
    \centering
    \caption{Comparison among the Static, Dynamic and Evolutionary strategies of GRASP with PR.}
\label{tab:StaticDynamicEvol}
\begin{tabular}{l c c c c c c c c c}
    \toprule
    & \multicolumn{4}{c}{$60$ s. } & &  \multicolumn{4}{c}{$600$ s. } \\
      \cline{2-5}  \cline{7-10} 
    \bf Procedure & \multicolumn{1}{c}{$\# Best$} & \multicolumn{1}{c}{$\% Dev$} & \multicolumn{1}{c}{$\# Best_k$} & \multicolumn{1}{c}{$\% Dev_k$}  & & \multicolumn{1}{c}{$\# Best$} & \multicolumn{1}{c}{$\% Dev$} & \multicolumn{1}{c}{$\# Best_k$} & \multicolumn{1}{c}{$\% Dev_k$} \\
   \midrule
   &  \multicolumn{9}{c}{LOP}\\
   \cline{2-10} 
    {\tt Static} &  61& 0.001& 60& 0.002& & 66& 0.001& 66& 0.001\\
    {\tt Dynamic}  & 53& 0.002& 53& 0.002& & 57& 0.002& 57& 0.002\\
    {\tt Evolutionary}  & 76& 0.000& 60& 0.001& & 78& 0.000& 70& 0.000\\
    \midrule
    
    & \multicolumn{9}{c}{MAX-CUT}\\
    \cline{2-10} 
      {\tt Static} &  53& 0.091& 17& 0.538& & 59& 0.076& 31& 0.300\\
    {\tt Dynamic}  & 13& 0.450& 14& 0.894& & 26& 0.310& 25& 0.532\\
    {\tt Evolutionary}  & 55& 0.096& 20& 0.543& & 43& 0.136& 30& 0.359\\
    \bottomrule
\end{tabular}

\end{table}

\begin{figure}[ht]
\centering
\includegraphics[width=0.8\textwidth]{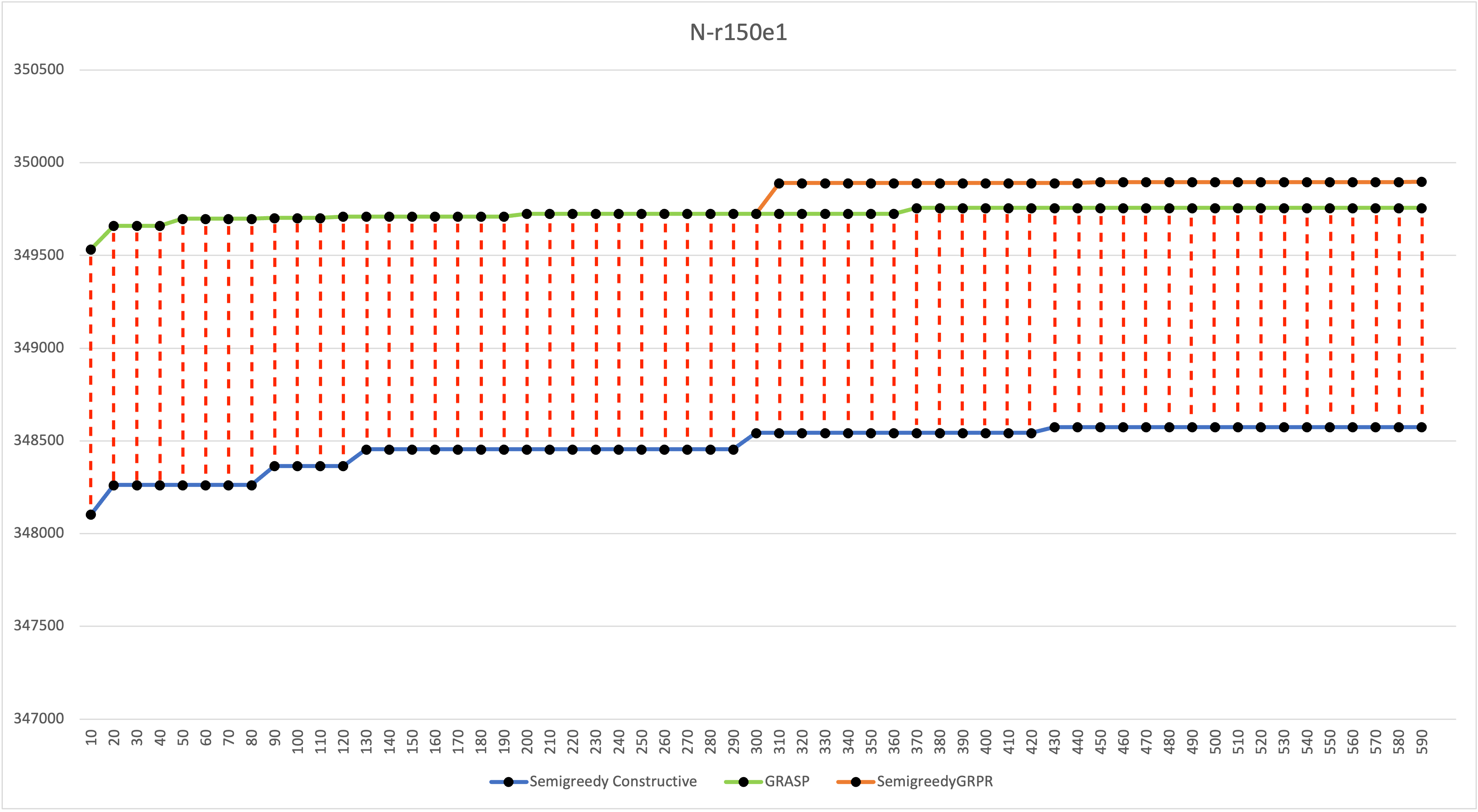}
\caption{Search profile for a representative instance of LOP.}
\label{fig:profile_lop}
\end{figure}

\begin{figure}[ht]
\centering
\includegraphics[width=0.8\textwidth]{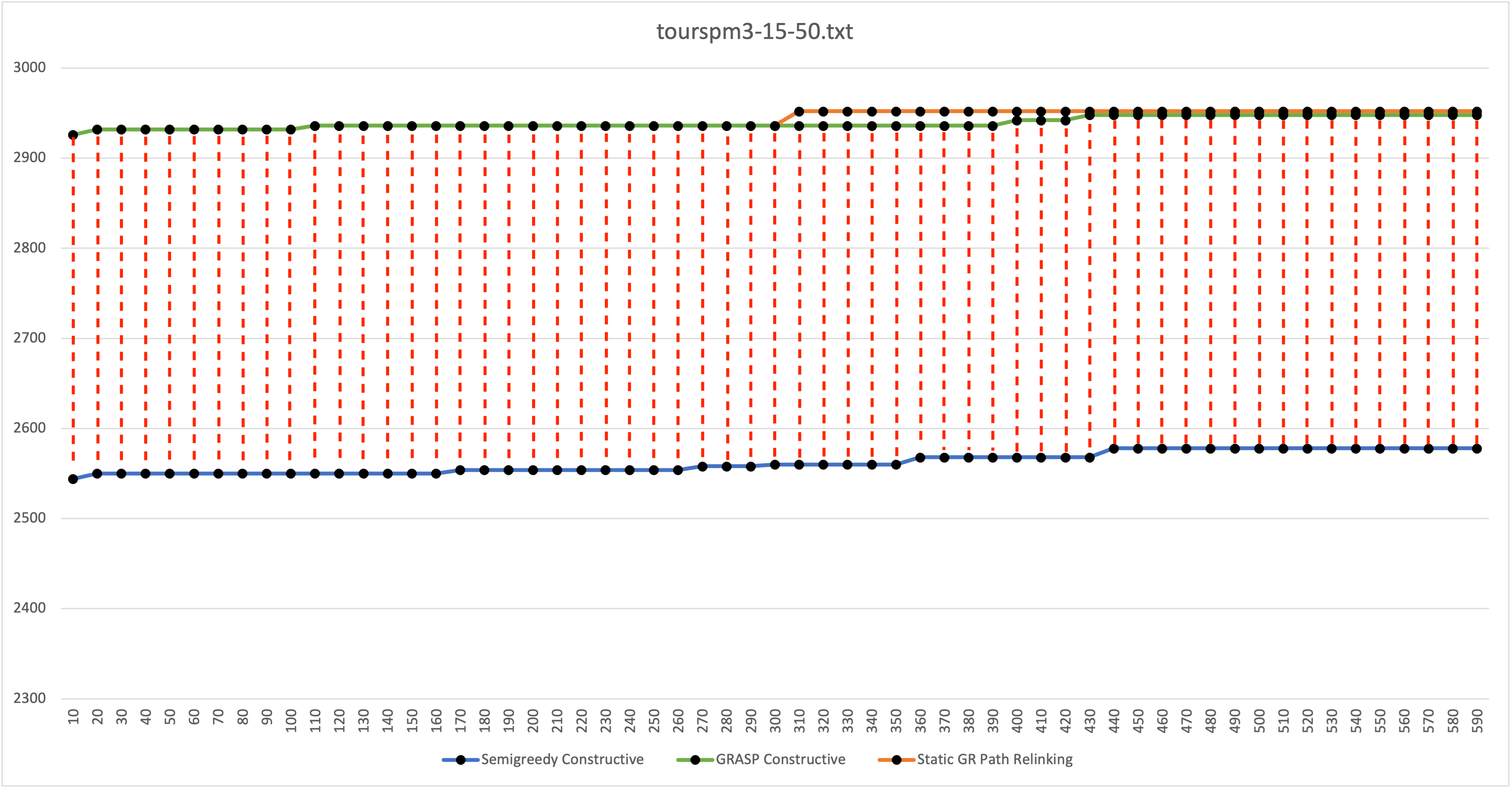}
\caption{Search profile for a representative instance of MAX-CUT.}
\label{fig:profile_maxcut}
\end{figure}

\underline{Experiment 7}. For our last experiment, we select one representative instance of LOP and one of MAX-CUT. We run three cases for 600 seconds: 1) construction method, 2) GRASP (i.e., construction and improvement), 3) GRASP for 300 seconds and then PR. The goal is to identify the contribution of the three main search element, namely, construction, improvement and PR. The blue line in Figures \ref{fig:profile_lop} and  \ref{fig:profile_maxcut} shows the objective function value of the incumbent solution throughout the run. The green line shows the objective function value of the incumbent solution for the GRASP run. The gap between these two lines (shown by the vertical red dashed lines) indicate the contribution of the improvement phase. The orange line shows the additional boost provided by PR at the 300-second mark. This boost is more evident in the LOP instance (Figure \ref{fig:profile_lop}) than in the MAX-CUT instance (Figure \ref{fig:profile_maxcut}). The experiment reveals a solution strategy in which there might be a  point in the search where the best course of action is to stop performing GRASP iterations and move to the PR phase.


\section{Conclusions}

We have provided a comprehensive review of the GRASP with Path Relinking metaheuristic and its significant impact on the field of combinatorial optimization over the last two decades. We have highlighted the main contributions of this algorithm including static, dynamic, and evolutionary strategies to hybridize GRASP and Path Relinking. Specifically, we  not only describe PR, but also how to construct and use a set of elite solutions in hybridizations with other metaheuristics, such as GRASP, that generate many local optimal solutions. In fact, the elite set serves as a long-term memory mechanism for GRASP in which PR is used for intensification. It is not our intent to be exhaustive in mentioning all contributions of GRASP with PR.  We limit our narrative to methodological contributions.

In our historical review, we observe that the first decade (2000-2010) was fruitful in terms of proposing new strategies and variants to hybridize GRASP and PR. In Section~\ref{sec:AdvancedPR} we focus on strategies and variants related to PR while in Section~\ref{sec:GRASP_PR} on those related to GRASP with PR. 
We discussed the different variations and extensions of the methodology proposed over these years, such as the pool management strategies for the elite set solutions generated with GRASP. 
In second decade (2010-2020) there was a consolidation of ideas where those strategies and variants were applied to solve new challenges.

We begin our experimental analysis in Section~\ref{sec:GraspExp} with a multi-start semi-greedy algorithm, and then add to it local search (resulting in a GRASP) and Path Relinking (resulting in a GRASP with PR). While there is no clear winner between the first and best improvement strategies in local search, in line with previous results in the literature, it is clear there is a benefit in adding local search to a semi-greedy heuristic.

In the second part of our experimental analysis in Section~\ref{sec:PRexp}, we compare forward, backward, and mixed strategies on a static GRASP with PR so that the affect of the different strategies on PR can be isolated. Our conclusion is that the best strategy is problem dependent where the mixed strategy performs better for the LOP while the forward strategy does better on the MAX-CUT.
We also studied the selective application of local search in intermediate solutions in the path of path relinking.
In terms of solution quality there appears to be little difference between three strategies (All, Every 5, and Only Best).  Since the Only Best strategy is the most efficient, this is recommended. 

Section~\ref{sec:GRASP-PRExp} compares two ways of selecting moves during the relinking process; a deterministic one called Greedy Path Relinking {\tt GPR}, and a randomized one called Greedy Randomized Path Relinking,  {\tt GRPR}. The results of this experiment are inconclusive in terms of an overall strategy for both problem types: while {\tt GRPR} is the better strategy for the LOP, {\tt GPR} is the better strategy for the MAX-CUT. In our final experiment analyzing the different elements and strategies of GRASP with PR, we consider different pool-management strategies to handle the elite set, namely Static, Dynamic, and Evolutionary, concluding that the Static and Evolutionary perform better than the Dynamic.  This seems to reinforce the merit of working with elite sets that include solutions of the highest quality.

Our analysis of various configurations and extensions of GRASP with Path Relinking indicate that their effectiveness varies across applications and context. We therefore recommend that analysts wishing to apply this technology investigate several of the strategies and elements that we have described here.  This could require the use of an automated algorithm configuration tool such as iterated racing \citep{lopez2016irace}.

Looking ahead, we anticipate that GRASP with Path Relinking will continue to play an important role in solving challenging optimization problems, and we expect to see further advances and refinements in its implementation and application. We hope that this review will inspire further research in this area and contribute to the ongoing development and improvement of this optimization technique.

\section{Acknowledgement}

This research has been partially supported by the Ministerio de Ciencia e Innovación of Spain (Grant Ref. PID2021-125709OB-C21 and PID2021-125709OA-C22) funded by MCIN/AEI/10.13039/ 501100011033/FEDER, UE. It has been also supported by the Generalitat Valenciana (CIAICO 2021/224).

\bibliographystyle{plainnat}
\bibliography{grasppr}

\pagebreak
\appendix
\section{Supplementary data} \label{append}

\begin{table}[htb]
    \centering
    \caption{Dynamic Path Relinking strategies for search direction.}
\label{tab:DynPathRelinking}
\begin{tabular}{l c c c c c c c c c}
    \toprule
    & \multicolumn{4}{c}{$60$ s. } & &  \multicolumn{4}{c}{$600$ s. } \\
      \cline{2-5}  \cline{7-10} 
    \bf Procedure & \multicolumn{1}{c}{$\# Best$} & \multicolumn{1}{c}{$\% Dev$} & \multicolumn{1}{c}{$\# Best_k$} & \multicolumn{1}{c}{$\% Dev_k$}  & & \multicolumn{1}{c}{$\# Best$} & \multicolumn{1}{c}{$\% Dev$} & \multicolumn{1}{c}{$\# Best_k$} & \multicolumn{1}{c}{$\% Dev_k$} \\
   \midrule
   &  \multicolumn{9}{c}{LOP}\\
   \cline{2-10} 
    {\tt DyPR Forward} &  66& 0.001& 49& 0.008& & 62& 0.001& 50& 0.005\\
    {\tt DyPR Backward}  & 54& 0.004& 48& 0.010& & 58& 0.003& 50& 0.007\\
    {\tt DyPR Mixed}  & 56& 0.008& 49& 0.015& & 58& 0.003& 51& 0.007\\
    \midrule
    
    & \multicolumn{9}{c}{MaxCut}\\
    \cline{2-10} 
    {\tt DyPR Forward} &  43& 0.208& 14& 1.221& & 41& 0.171&  15& 0.807\\
    {\tt DyPR Backward}  &  43& 0.182& 14& 1.194& & 53& 0.146& 18& 0.782\\
    {\tt DyPR Mixed}  &  46& 0.163& 13& 1.176& & 43& 0.227& 13& 0.863\\
    \bottomrule
\end{tabular}

\end{table}

\begin{table}[htb]
    \centering
    \caption{Evolutionary Path Relinking strategies for search direction.}
\label{tab:EvPathRelinking}
\begin{tabular}{l c c c c c c c c c}
    \toprule
    & \multicolumn{4}{c}{$60$ s. } & &  \multicolumn{4}{c}{$600$ s. } \\
      \cline{2-5}  \cline{7-10} 
    \bf Procedure & \multicolumn{1}{c}{$\# Best$} & \multicolumn{1}{c}{$\% Dev$} & \multicolumn{1}{c}{$\# Best_k$} & \multicolumn{1}{c}{$\% Dev_k$}  & & \multicolumn{1}{c}{$\# Best$} & \multicolumn{1}{c}{$\% Dev$} & \multicolumn{1}{c}{$\# Best_k$} & \multicolumn{1}{c}{$\% Dev_k$} \\
   \midrule
   &  \multicolumn{9}{c}{LOP}\\
   \cline{2-10} 
    {\tt EvPR Forward} &  65& 0.001& 49& 0.007& & 62& 0.001& 51& 0.005\\
    {\tt EvPR Backward}  & 59& 0.004& 48& 0.010& & 59& 0.001& 51& 0.005\\
    {\tt EvPR Mixed}  & 53& 0.009& 49& 0.016& & 59& 0.002& 52& 0.006\\
    \midrule
    
    & \multicolumn{9}{c}{MaxCut}\\
    \cline{2-10} 
    {\tt EvPR Forward} &  51& 0.109& 16& 0.870& & 45& 0.138&  26& 0.474\\
    {\tt EvPR Backward}  &  38& 0.219& 15& 0.978& & 51& 0.134& 27& 0.470\\
    {\tt EvPR Mixed}  &  38& 0.235& 13& 0.994& & 34& 0.263& 16& 0.599\\
    \bottomrule
\end{tabular}

\end{table}

\end{document}